\documentclass[a4paper,oneside, 12 pt, reqno]{amsart}
\usepackage{enumerate}
\usepackage{amsmath,amsthm,amscd,amsfonts,amssymb,mathrsfs}
\usepackage{latexsym,amsmath}
\usepackage{mathrsfs}
\usepackage[maxbibnames=99, backend=bibtex]{biblatex}
\usepackage{xcolor}
\usepackage{soul}
\usepackage{comment}
\makeatletter
\@namedef{subjclassname@2010}{%
\textup{2010} Mathematics Subject Classification}
\makeatother

\newcommand{\ncom}{\newcommand}
\ncom{\beqn}{\begin{eqnarray*}}
	\ncom{\eeqn}{\end{eqnarray*}}
\ncom{\beq}{\begin{eqnarray}}
	\ncom{\eeq}{\end{eqnarray}}
\ncom{\cal}{\mathcal}
\ncom{\eop}{\hfill{{\rule{2.5mm}{2.5mm}}}}
\ncom{\eoe}{\hfill{{\rule{1.5mm}{1.5mm}}}}
\ncom{\eof}{\hfill{{\rule{1.5mm}{1.5mm}}}}
\ncom{\hone}{\mbox{\hspace{1em}}}
\ncom{\htwo}{\mbox{\hspace{2em}}}
\ncom{\hthree}{\mbox{\hspace{3em}}}
\ncom{\hfour}{\mbox{\hspace{4em}}}
\ncom{\hsev}{\mbox{\hspace{7em}}}
\ncom{\vone}{\vskip 2ex}
\ncom{\vtwo}{\vskip 4ex}
\ncom{\vonee}{\vskip 1.5ex}
\ncom{\vthree}{\vskip 6ex}
\ncom{\vfour}{\vspace*{8ex}}
\ncom{\norm}{\|\;\;\|}
\ncom{\integ}[4]{\int_{#1}^{#2}\,{#3}\,d{#4}}
\ncom{\inp}[2]{\langle{#1},\,{#2} \rangle}
\ncom{\Inp}[2]{\Langle{#1},\,{#2} \Langle}
\ncom{\vspan}[1]{{{\rm\,span}\#1 \}}}
\ncom{\dm}[1]{\displaystyle {#1}}

\newtheorem{theorem}{\bf Theorem}[section]
\newtheorem{lemma}[theorem]{\bf Lemma}%[section]

\newtheorem{question}[theorem]{\bf Question}

\newtheoremstyle
{remarkstyle}
{}
{11pt}
{}
{}
{\bfseries}
{:}
{     }
{\thmname{#1} \thmnumber{#2} }

\theoremstyle{remarkstyle}

\newtheorem{cor}[theorem]{Corollary}

\begin{document}

	\title[CDSP for cyclic, analytic $2$-isometry]{\bf A solution to the Cauchy dual subnormality problem for a cyclic analytic $2$-isometry with defect operator of rank two}
	
	\author[a1]{Mandar Khasnis}
	\address{Department of Mathematics, Smt. CHM College, Ulhasnagar, 421003}
	\email{mkhasnis.official@gmail.com}
	
	\author[a2]{Geetanjali Phatak}
	\address{Department of Mathematics, S. P. College, Pune, 411030}
	\email{gmphatak19@gmail.com}
	
	\author[a3]{V. M. Sholapurkar*}
	\address{Bhaskaracharya Pratishthan, Pune, 411004}
	\email{*Corresponding Author: vmshola@gmail.com}
	
	\date{}

	\begin{abstract}
		
		The Cauchy dual subnormality problem (for short, CDSP) asks whether the Cauchy dual of a $2$-isometry is subnormal. In this article, we prove that if $\mu$ is a sum of unit point mass measures at two non-antipodal points on the unit circle, then the Cauchy dual $M_z'$ of the multiplication operator $M_z$ on the Dirichlet-type space $D(\mu)$ is not subnormal. If the points are antipodal then the subnormality of the said operator has been already established in the literature. Thus, we have a complete solution to CDSP in this case.
	\end{abstract}
	\subjclass[2020]{47B38, 47B20, 47B32}
	\keywords{Cauchy dual; Subnormal Operator; de Branges–Rovnyak space; Dirichlet-type spaces\\ *Corresponding Author: V. M. sholapurkar. Email: vmshola@gmail.com}
	\maketitle
	%\maketitle
	%\tableofcontents
	%\setcounter{tocdepth}{1}

	\section{Introduction}
	In this article, we present a solution to the Cauchy dual subnormality problem (for short, CDSP) for a family of cyclic, analytic $2$-isometries $T$ whose defect operator $T^*T-I$ is of rank two. \\ 
	The notion of Cauchy dual of an operator was introduced by Shimorin in \cite{Sh}. However, the genesis of the concept can be traced in the following classical result. 
	\begin{theorem}\label{t1}
		The reciprocal of a Bernstein function is completely monotone. 
	\end{theorem}
	A proof of Theorem \ref{t1} and the wealth of information on various special classes of functions related to Bernstein functions can be found in \cite{bernstein_ssv2009}. The interplay between the classes of functions described in Theorem \ref{t1} gets mirrored in the classes of operators that can be associated with these classes of functions. In fact, the {\it Cauchy dual of an operator} is the transform that allows one to carry this function theoretic correspondence to an operator theoretic relationship. \\
	An excellent illustration of such a duality is provided by {\it contractive subnormal operators} and {\it completely hyperexpansive operators}.  The interplay between these two classes of operators can be thought of as an operator-theoretic manifestation of two special classes of functions, viz. completely monotone functions and completely alternating functions. This fruitful synthesis of harmonic analysis with operator theory has been extensively explored in the works of Athavale et al. \cite{athavale2003, athavale_ranjekar2002_1, athavale_ranjekar2002_2, athavale_vms1999,  sholapurkar2000}, as well as in the works of  Jab{\l}o{\'n}ski and Stochel \cite{jablonski2006, jablonski2002, jablonski2003, jablonski2001, jablonski2004}. It turns out that the notion of Cauchy dual fits naturally in these circumstances and helps in understanding the subtle interconnections.
	
	The said theme begins with an observation: 
%	If $T:\{\alpha_n\}$ is a weighted shift operator, then the Cauchy dual $T'$ is the weighted shift with weight sequence  $\Big\{\frac{1}{\alpha_n}\Big\}.$ Combining this fact with a result in \cite[Chapter 4, Proposition 6.10]{bcr1984}  enabled Athavale to prove that 
	The Cauchy dual of a completely hyperexpansive weighted shift is a contractive subnormal weighted shift \cite{athavale1996}. In the light of this observation, we naturally have the following question:
	\begin{question}\label{q}
		{\bf Is the Cauchy dual of a completely hyperexpansive operator a subnormal contraction?}
	\end{question}
	This question appears in \cite[Question 2.11]{Ch}. 
	%Some special cases of the question have been dealt with 
	The question has been extensively studied in the literature. The following results indicate that the evidence is in the positive direction: 
	\begin{enumerate}
		\item The Cauchy dual of a completely hyperexpansive weighted shift is a subnormal weighted shift \cite{athavale1996}.
		\item The Cauchy dual of a concave operator is a hyponormal contraction \cite{shimorin2002}.
		\item The Cauchy dual of a $\Delta_T$-regular 2-isometry is subnormal  \cite{badea2019}.
		\item The Cauchy dual of a 2-isometric operator-valued weighted shift is subnormal \cite{acjs2019}.
		\item The Cauchy dual of an expansive operator of class Q is  subnormal \cite{CJJS2021}.
	\end{enumerate}
	In view of the above results, the following counterexamples are instructive.
	\begin{enumerate}
		\item There exist weighted shifts on directed trees whose Cauchy dual is not subnormal \cite{acjs2019}.
		\item There exist a class of cyclic $2$-isometric composition operators whose Cauchy dual is not subnormal \cite{sc2020}.
	\end{enumerate}
	It is well known that the class of $2$-isometric operators is a subclass of the class of completely hyperexpansive operators. The results stated above suggest that even the following special case of the  Question \ref{q}  also seems to be out of sight:   
	\begin{question}\label{q1}
		Is the Cauchy dual of a cyclic $2$-isometric operator a subnormal operator? 
	\end{question}
	
	% In this circle of ideas, we record the following results dealing with the variants of the CDSP.                
	%\begin{theorem}(\cite{Ch}, theorem 2.9)
	
	%The operator Cauchy dual of a 2-hyperexpansion is a hyponormal contraction. 	
	% \end{theorem}

The Question \ref{q1} can be reformulated in terms of the classification of finite, positive, Borel measures on the unit circle. Recall that for a finite positive Borel measure $\mu$ on the unit circle $\mathbb{T}$, the Dirichlet-type space $D(\mu)$ is defined by
\[D(\mu):=\Bigg\{f\in \text{Hol}(\mathbb{D}) : \displaystyle\int_{\mathbb{D}}|f'(z)|^2P_\mu(z)dA(z)<\infty \Bigg\}.\]
Here, $P_\mu(z)$ is the Poisson integral for measure $\mu$ given by $\displaystyle\int_{\mathbb{T}}\displaystyle\frac{1-|z|^2}{|z-\zeta|^2}d\mu(\zeta)$ and $dA$ denotes the normalized Lebesgue area measure on the open unit disc $\mathbb{D}$. %and Hol$(\mathbb{D})$ denotes the set of all holomorphic functions on $\mathbb{D}$. 
The space $D(\mu)$ is a reproducing kernel Hilbert space and multiplication by the coordinate function $z$, denoted by $M_z$, turns out to be a bounded linear operator (see \cite{aleman1993,richter1991}) on $D(\mu)$. For an elaborate discussion on Dirichlet-type spaces, the reader may also consult \cite{dirichlet_sp2014}. S. Richter described a model for a cyclic, analytic $2$-isometry \cite{richter1991}. Indeed, a cyclic, analytic $2$-isometry is unitarily equivalent to a multiplication operator $M_z$ on a Dirichlet-type space $D(\mu)$ for some finite positive Borel measure $\mu$ on unit circle $\mathbb{T}$. This model for cyclic analytic 2-isometries and a decomposition theorem of 2-isometries allow us to revisit the the following question that appeared in  \cite{cgr2022}:  
\begin{question}\label{q2}
	Classify all finite positive Borel measures on the unit circle for which the Cauchy dual $M_z^{\prime}$ of the multiplication operator $M_z$ on the Dirichlet-type space $D(\mu)$ is subnormal.
\end{question}
It was observed in \cite{cgr2022} that the rank of the defect operator is infinite for all the counterexamples to CDSP that appeared in the literature. Therefore,  the following particular case of Question \ref{q1} was posed in \cite{cgr2022}: 
\begin{question}\label{q3} 
	Is it true that for all cyclic, analytic $2$-isometries $T$ for which $T^*T-I$ is of finite rank, the Cauchy dual of $T$ is subnormal? 
\end{question}
In fact, in the same article \cite{cgr2022}, the authors constructed an example of a cyclic, analytic $2$-isometry $T$ for which $T^*T-I$ is of rank two whose Cauchy dual is subnormal. The example certainly strengthened the evidence for the affirmative answer to Question \ref{q3}. However, recently an example of a cyclic, analytic $2$-isometry $T$ for which $T^*T-I$ is of rank two whose Cauchy dual is not subnormal has been found in \cite{mkvms2025}.  \\ 
At this stage, the reader is referred to the fact \cite[Theorem 6.1]{cgr2022} which states that for a cyclic, analytic $2$-isometry $T$ in $\cal B(H)$, the rank of $T^*T-I$ is finite if and only if there exists a {\it finitely supported} measure $\mu$ on the unit circle $\mathbb{T}$  such that $T$ is unitarily equivalent to the operator $M_z$ on $D(\mu)$. In fact, if the rank of $T^*T-I$ is a positive integer $k,$ then the corresponding measure $\mu$ is supported at exactly $k$ points on the unit circle $\mathbb{T}.$
In view of this result, a construction of a cyclic, analytic $2$-isometry with defect operator of finite rank  reduces to choosing finitely many points on the unit circle and looking at the multiplication operator $M_z$ on the Dirichlet-type space $D(\mu)$, where $\mu$ is supported at the chosen points. 
Further, in \cite[Theorem 6.4]{cgr2022} it is claimed that this space $D(\mu)$ coincides with the de Branges-Rovnyak space $H(B)$ for some vector-valued rational function $B$ with equality of norms. \\
We now quote two examples, which led us to guess the solution of the problem presented in this article.
\begin{itemize}
	\item An example of a cyclic, analytic $2$-isometry $T$ such that $T^*T-I$ is of rank two and whose Cauchy dual is {\it subnormal} has been constructed in \cite [Example 7.2] {cgr2022} by choosing a measure on the unit circle supported at $\{1,-1\}$.
	\item An example of a cyclic, analytic $2$-isometry $T$ such that $T^*T-I$ is of rank two and whose Cauchy dual is {\it not subnormal} has been constructed in \cite{mkvms2025} by choosing a measure on the unit circle supported at $\{1,i\}$.
\end{itemize}

In view of the examples cited above, the Question \ref{q2} can be narrowed down to the following question: 
\begin{question}
	Classify all finitely supported positive Borel measures on the unit circle for which the Cauchy dual $M_z^{\prime}$ of the multiplication operator $M_z$ on the Dirichlet-type space $D(\mu)$ is subnormal and $M_z^*M_z-I$ is of rank two.
\end{question}
The present article is devoted to a partial solution of the above question. The main theorem proved in this paper is as follows. 
\begin{theorem}[Main Theorem]\label{main}
	Let $\zeta_1, \zeta_2$ be two distinct points on the unit circle. Let $\mu$ be a positive Borel measure of the form $\mu= \delta_{\zeta_1}+\delta_{\zeta_2}$ defined on the unit circle $\mathbb{T}$. If $\zeta_1, \zeta_2$ are non-antipodal points, then the Cauchy dual of $M_z$ on $D(\mu)$ is not subnormal.
\end{theorem}
The proof of the main theorem is presented in section \ref{sec3}. Note that in view of \cite[Proposition 7.1]{cgr2022}, it is sufficient to prove the theorem in the case where the measure $\mu$ is supported on $\{1, \xi\} $ for any point $\xi~(\xi\neq -1) $ on the unit circle. In the process, we shall obtain the reproducing kernel for the space $D(\mu)$ where $\mu$ is supported on $\{1,\xi\}$. \\
Theorem \ref{main} coupled with \cite [Theorem 2.4] {cgr2022} yields the following characterization in this case.
\begin{theorem}
	Let $\zeta_1, \zeta_2$ be two distinct points on the unit circle. Let $\mu$ be a positive Borel measure of the form $\mu= \delta_{\zeta_1}+\delta_{\zeta_2}$ defined on the unit circle $\mathbb{T}$. Then the Cauchy dual of $M_z$ on $D(\mu)$ is subnormal if and only if the points $\zeta_1$ and $\zeta_2$ are antipodal.  
\end{theorem}

\section{Preliminaries}
Let $\mathbb{R}, \mathbb{C}, \mathbb{D}$ and $\mathbb{T}$, respectively, denote the set of real numbers, the set of complex numbers, the open unit disc and the unit circle.

The notion of a reproducing kernel Hilbert space (RKHS) is well known and has been extensively studied, explored and crucially used in the theory of operators on Hilbert spaces. We quickly recall the definition of RKHS for a ready reference. 
%{\bf Reproducing Kernel Hilbert Space}:\\

Let $S$ be a non-empty set and $\cal F$$(S,\mathbb{C})$ denote set of all functions from $S$ to $\mathbb{C}$. We say that $\cal H\subseteq \cal F$$(S,\mathbb{C})$ is a {\it \bf Reproducing Kernel Hilbert Space}, briefly {\it RKHS}, if
\begin{enumerate}
	\item[i.] $\cal H$ is a vector subspace of $\cal F$$(S,\mathbb{C})$
	\item[ii.] $\cal H$ is equipped with an inner product $\langle , \rangle$ with respect to which $\cal H$ is a Hilbert space
	\item[iii.] For each $x\in S$, the linear evaluation functional $E_x:\cal H$$\to \mathbb{C}$ defined by $E_x(f)=f(x)$, is bounded.
\end{enumerate}
If $\cal H$ is an RKHS on $S$ then for each $x\in S$, by Riesz representation theorem, there exists unique vector $k_x\in \cal H$ such that, for each $f\in \cal H$,\\ $f(x)=E_x(f)=\langle f,k_x\rangle$.
\\The function $K:S\times S \to \mathbb{C}$ defined by $K(x,y)=k_y(x)$ is called as the reproducing kernel for $\cal H$. A classical reference for reproducing kernel Hilbert spaces is \cite{paulsen2016RKHS}.

We now recall the definitions of three special types of RKHS which are relevant in the present context. 
\begin{enumerate}
	
	\item {\bf Hardy Space:} Let Hol$(\mathbb{D})$ denote the set of all holomorphic functions on $\mathbb{D}$. The Hardy space $H^2$ is defined as
	\[H^2=\Bigg\{f(z)=\displaystyle\sum_{n=0}^{\infty}a_nz^n\in \text{Hol}(\mathbb{D}) : \displaystyle\sum_{n=0}^{\infty}|a_n|^2<\infty\Bigg\}.\]
	It is a well known fact that $H^2$ is a reproducing kernel Hilbert space with the inner product defined as $\langle f,g \rangle :=\displaystyle\sum_{n=0}^{\infty}a_n\overline{b_n}$, if $f(z)=\displaystyle\sum_{n=0}^{\infty}a_nz^n$ , $g(z)=\displaystyle\sum_{n=0}^{\infty}b_nz^n$ and a kernel (namely {\it Szego Kernel}) function $K_\lambda(z)=\displaystyle\frac{1}{1-\overline{\lambda}z}$. The interested reader may consult \cite{douglas2012}.\\
	
	\item {\bf Dirichlet-type spaces:}  
	Note that Dirichlet-type spaces have been described in the previous section.
	D. Sarason \cite{sarason97} initiated the study of Dirichlet-type spaces by way of identifying such a space with a de Branges-Rovnyak space. This association has been further strengthened in recent years in the works of  \cite{cgr2022, cgr2010, kellay2015}. The identification of a Dirichlet-type space as a de Branges-Rovnyak space allows one to compute the reproducing kernel for the Dirichlet-type space. \\ Here, we include a brief description of de Branges-Rovnyak space for a ready reference. 
	
	\item {\bf de Branges-Rovnyak spaces:}
	\\For complex, separable Hilbert spaces $\cal U$, $\cal V$, let $\cal B(U,V)$ denote the Banach space of all bounded linear transformations from $\cal U$ to $\cal V$. The {\bf Schur class} $S(\cal U,V)$ is given by
	\[S({\cal U,V})=\Big\{B:\mathbb{D} \to {\cal B(U,V)} : B \text{ is holomorphic, } \displaystyle\sup_{z\in \mathbb{D}}\|B(z)\|_{B({\cal U,V})}\leq 1\Big\}.\]
	Observe that when ${\cal U}=\mathbb{C}=\cal V$ then the Schur class is nothing but the closed unit ball of $H^\infty(\mathbb{D})$, the set of all bounded holomorphic functions on $\mathbb{D}$. For any $B\in S(\cal U,V)$, the de-Branges-Rovnyak space, $H(B)$ is the reproducing kernel Hilbert space associated with the $\cal B(V)$-valued semidefinite kernel given by
	\[K_B(z,w)=\displaystyle\frac{I_{\cal V}-B(z)B(w)^*}{1-z\overline{w}},~~z,w\in \mathbb{D}.\]
	For equivalent formulations of de Branges-Rovnyak spaces, the reader is referred to an excellent article by J. Ball \cite{Ball2015}. 
	The kernel $K_B$ is normalised if $K_B(z,0)=I_{\cal V}$ for every $z\in \mathbb{D}$. This is equivalent to the condition $B(0)=0$. Further, when ${\cal U}=\mathbb{C}=\cal V$, we denote $H(B)$ by $H(b)$, the classical de Branges-Rovnyak space (Refer \cite{fricain_mashreghi_2016vol2} for the basic theory of the classical de Branges-Rovnyak spaces). 
\end{enumerate}
We now record a few definitions which are required in the sequel.\\
Let $\cal H$ be a complex, infinite dimensional, separable Hilbert space and $\cal B(\cal H)$ denotes the $C^*$-algebra of bounded linear operators on $\cal H$. 

Let $T\in {\cal B(H)}$.
We say that $T$ is {\it cyclic} if there exists a vector $f\in \cal H$ such that ${\cal H}=\bigvee \big\{T^nf : n\geq 0\big\}$ (vector $f$ is known as {\it cyclic vector}). An operator $T$ is said to be {\it analytic} if $\cap_{n\geq 0}T^n{\cal H}=\{0\}$. The  {\it Cauchy dual} $T'$ of a left invertible operator $T\in \cal B(\cal H)$ is defined as $T'=T(T^*T)^{-1}$. Following Agler \cite{agler1995_1}, an operator $T\in \cal B(\cal H)$ is said to be a {\it $2$-isometry} if \[I-2T^*T+T^{*2}T^2=0.\]
An operator $T\in \cal B(H)$ is said to be {\bf subnormal } if there exist a Hilbert space $\cal K$ and an operator $S\in \cal B(K)$ such that $\cal H \subseteq K$, $S$ is normal and $S|_{\cal H}=T$. Readers are encouraged to refer \cite{conway1991} for a detailed study of subnormal operators. Agler proved in \cite{agler1985} that $T\in \cal B(H)$ is a subnormal contraction ($T$ is subnormal and $\|T\|\leq 1$) if and only if 
\[B_n(T)\equiv \sum_{k=0}^{n}(-1)^k \binom{n}{k}T^{*k}T^k\geq 0 \text{ for all integers }n\geq 0.\]
Following \cite{athavale1996}, an operator $T\in B(\cal H)$ is said to be {\it completely hyperexpansive} if $B_n(T)\leq 0~~ \text{for all integers }n\geq 1.$ 

We now record the results used in the proof of main theorem for a ready reference.   
\subsection{Reproducing kernel of $D(\mu)$}\label{sec2}
The algorithm for obtaining reproducing kernel of $D(\mu)$ turns out to be a fruitful computational tool for checking the subnormality of the Cauchy dual of $M_z$ on $D(\mu)$. For $\xi \in \mathbb{T}~(\xi\neq \pm 1)$, we take up the task of computing the reproducing kernel for the space $D(\mu)$ where $\mu$ is the regular, positive, Borel measure on the unit circle $\mathbb{T}$ of the form $\delta_{1}+\delta_{\xi},$ where $\delta_\lambda$ denotes the unit point mass measure at a point $\lambda\in \mathbb{T}.$ The technique used in obtaining the reproducing kernel heavily relies on the work of Costara \cite{costara2016}.  Thus, we briefly outline the work carried out in \cite{costara2016} and state the results that are used in the sequel.
Given $n\in \mathbb{N}$, let $c_1,c_2,\dots,c_n$ be positive real numbers and $\zeta_1,\zeta_2,\dots,\zeta_n$ be distinct points on the unit circle $\mathbb{T}$. Consider the positive Borel measure  $\displaystyle \mu=\sum_{j=1}^{n}c_j\delta_{\zeta_j}$ on $\mathbb{T}$.\\ %where $\delta_{\zeta_j}$ denotes the Dirac delta measure supported at $\zeta_j$.\\
As an application of Riesz-Fejer theorem, it is observed in \cite{costara2016} that there exist $\alpha_1,\alpha_2,\dots,\alpha_n \in \mathbb{C}\setminus\overline{\mathbb{D}}$ and $d>0$ such that the following equation holds:
\begin{equation}\label{eq3}
	\prod_{j=1}^{n}|z-\zeta_j|^2+\sum_{j=1}^{n}c_j\prod_{\substack{i=1\\i\neq j}}^{n}|z-\zeta_i|^2=d\prod_{j=1}^n |z-\alpha_j|^2,~~~z\in \mathbb{T}.
\end{equation}
The values of $\alpha_1,\alpha_2,\dots,\alpha_n$ as obtained in equation (\ref{eq3}) turn out to be important in the expression for reproducing kernel of $H^2.$ Following \cite[Theorem 3.1, Theorem 4.4]{costara2016},
\begin{enumerate}
	\item Define \begin{equation}\label{eq16}
		O_\mu(z)=\displaystyle\frac{p(z)}{q(z)}
	\end{equation}
	with $p(z)=\displaystyle\frac{e^{i\theta}}{\sqrt{d}}\displaystyle\prod_{j=1}^{n}(z-\zeta_j)$ ($\theta \in \mathbb{R}$ is chosen such that $O_\mu(0)>0$) and  $q(z)=\displaystyle\prod_{j=1}^{n}(z-\alpha_j)$.\\
	\noindent We observe that $O_\mu$ is a rational function of degree $n$ having poles outside the closed unit disk and simple zeros at $\zeta_j\in \mathbb{T}$, for each $j=1,\dots,n$. Then $\prod_{j=1}^{n}(z-\zeta_j)H^2=O_\mu H^2$ and hence $O_\mu H^2$ is a closed subspace of $D(\mu)$.
	The reproducing kernel for $O_\mu H^2$ is given by (Refer \cite[Theorem 3.1]{costara2016}) 
	\begin{equation}\label{eq4}
		\tilde{K_\mu}(z,\lambda)=\displaystyle \frac{O_\mu(z)}{1-\overline{\lambda}z}\overline{O_\mu(\lambda)}.
	\end{equation}
	\item The reproducing kernel for $(O_\mu H^2)^\perp$ as a subspace of $D(\mu)$ is given by (Refer \cite[Theorem 4.4]{costara2016})
	\begin{equation} \label{eq5}
		\hat{K_\mu}(z,\lambda)=\displaystyle\sum_{j=1}^{n}g_j(\lambda)f_j(z) ~~~~ z,\lambda \in \mathbb{D}
	\end{equation}\label{eq6}
	where \begin{equation}\label{eq11}
		f_j(z)=\displaystyle\frac{O_\mu(z)}{O_\mu '(\zeta_j)(z-\zeta_j)}
	\end{equation}
	and 
	\begin{equation}\label{eq7}
		\begin{pmatrix}
			\overline{g_1(\lambda)}\\
			\overline{g_2(\lambda)}\\
			\vdots \\
			\overline{g_n(\lambda)}
		\end{pmatrix}=\left(\left[\langle f_i,f_j\rangle_\mu\right]_{1\leq i,j\leq n}\right)^{-1}
		\begin{pmatrix}
			f_1(\lambda)\\
			f_2(\lambda)\\
			\vdots \\
			f_n(\lambda)
		\end{pmatrix}.
	\end{equation}	
	\item The Lemma \cite[Lemma 4.2, 4.3]{costara2016} provides formula for obtaining $\langle f_i,f_j\rangle_\mu$ for any $1\leq i,j\leq n$ as follows :
	\begin{equation}\label{eq8}
		\|f_i\|_\mu^2=\langle f_i,f_i\rangle_\mu=c_i\zeta_i f_i'(\zeta_i)
	\end{equation} and for $i\neq j$,
	\begin{equation}\label{eq9}
		\langle f_i,f_j\rangle_\mu= \displaystyle\frac{1}{O_\mu'(\zeta_i)\overline{O_\mu'(\zeta_j)}(1-\zeta_i\overline{\zeta_j})}.
	\end{equation}
\end{enumerate}

\noindent The implementation of the above algorithm by choosing  $\zeta_1=1$, $\zeta_2=\xi$ and $c_1=c_2=1$ yields the following theorem:
\begin{theorem}\label{thm1}
	Let $\mu$ be a positive Borel measure on unit circle $\mathbb{T}$ defined as $\mu=\delta_1+\delta_\xi$. The reproducing kernel for $D(\mu)$ is given by:
	\begin{align*}
		K(z,\lambda)&
		=\displaystyle\frac{b}{\overline{q(\lambda)}q(z)}\Bigg[\displaystyle\frac{(z-1)(z-\xi)(\bar \lambda -1)(\bar\lambda -\bar\xi)}{1-\bar\lambda z}\\ &+(a-1)\left[(z-\xi)(\bar\lambda -\bar \xi)+(z-1)(\bar\lambda-1) \right]\\ 
		&+\displaystyle\frac{(z-\xi)(\bar\lambda-1)}{\Delta(\xi-1)}+\frac{(z-1)(\bar\lambda - \bar\xi)}{\Delta(\bar\xi-1)}
		\Bigg], ~~z,\lambda\in \mathbb{D}
	\end{align*}
	where $a,b, \Delta\in \mathbb{R}$ are constants and $q(z)$ is a polynomial of degree $2$ whose roots are in $\mathbb{C}\setminus \overline{\mathbb{D}}$.
	%$\Delta=\det\left(\langle f_i,f_j\rangle_\mu\right)_{1\leq i,j\leq 2}$.
		%\begin{align*}
	%K(z,\lambda)& =\displaystyle\frac{b}{\overline{q(\lambda)}q(z)}\Bigg[\displaystyle\frac{(\bar \lambda -1)(\bar\lambda +i)(z-1)(z-i)}{1-\bar\lambda z}\\ &+(a-1)\left((\bar\lambda +i)(z-i)+(\bar\lambda -1)(z-1)\right)\\ &+\displaystyle\frac{\bar s(\bar \lambda -1)(z-i)}{-i(w-pi)^2}+\frac{s(\bar \lambda +i)(z-1)}{i(w+pi)^2}\Bigg] ~~~~~~~(\lambda, z\in \mathbb{D})
	%\end{align*} 
	%where $a,b,p,w,s\in \mathbb{R}$ are constants with $p^2+w^2=1$ and $q(z)$ is a polynomial of degree $2$ whose roots are in $\mathbb{C}\setminus \overline{\mathbb{D}}$.\\
\end{theorem}
We now proceed towards the proof of the main theorem.
\section{Proof of the Main Theorem}\label{sec3}
It is proved in \cite[Theorem 2.4]{cgr2022} that in the case when measure $\mu$ is supported at a single point, the  Cauchy dual of the operator $M_z$ on $D(\mu)$ is subnormal. Also, as mentioned earlier, if the measure is supported at $\{1,-1\}$, the Cauchy dual of $M_z$ on $D(\mu)$ is subnormal. Thus, we set $\xi=\cos \theta +i \sin \theta$ where $\theta \neq 0, \pi$. The proof of the main theorem is an application of the following theorem from \cite{cgr2022}.
\begin{theorem}\cite[Theorem 6.4]{cgr2022}\label{thm2}
	For positive scalars $c_1,c_2,\dots,c_k$ and distinct points $\zeta_1,\zeta_2,\dots,\zeta_k$ on the unit circle $\mathbb{T}$, consider the positive Borel measure $\mu=\sum_{j=1}^{k}c_j\delta_{\zeta_j}$ on $\mathbb{T}$, where $\delta_{\zeta_j}$ denotes the Dirac delta measure supported at $\zeta_j$. Let $X(z)=(z,z^2,\dots,z^k)^T$ and $\{e_j\}_{j=1}^k$ denote the standard basis of $\mathbb{C}^k$. Then there exist $\alpha_1, \alpha_2,\dots,\alpha_k\in \mathbb{C}\setminus\overline{\mathbb{D}}$ and a $k\times k$ upper triangular matrix $P$ such that the Dirichlet-type space $D(\mu)$ coincides with the de Branges-Rovnyak space $H(B)$ with equality of norms, where $B=\left(\frac{p_1}{q}, \frac{p_2}{q}, \dots, \frac{p_k}{q}\right)$ and 
	$p_j(z)=\left\langle PX(z),e_j\right\rangle~,~j=1,\dots,k,~~~q(z)=\prod_{j=1}^k (z-\alpha_j).$
	Moreover, $\alpha_1,\dots,\alpha_k$ are governed by
	\[\prod_{j=1}^{k}|z-\zeta_j|^2+\sum_{j=1}^{k}c_j\prod_{\substack{i=1\\i\neq j}}^{k}|z-\zeta_i|^2=d\prod_{j=1}^k |z-\alpha_j|^2,~~~z\in \mathbb{T}\]
	for some $d>0$.
\end{theorem}
In addition, we record the following two results from \cite{cgr2022} that are used in the sequel.
\begin{theorem}\cite[Theorem 2.1]{cgr2022}\label{thm3}
	Let $B=(b_1,\dots,b_k)\in {\cal S}(\mathbb{C}^k,\mathbb{C})$ be such that $B(0)=0$ where
	\[b_j(z)=\displaystyle\frac{p_j(z)}{\prod_{j=1}^k (z-\alpha_j)}\]
	for polynomials $p_j$ of degree at most $k$ and distinct numbers $\alpha_1,\dots,\alpha_k \in \mathbb{C}\setminus\overline{\mathbb{D}}$. For $r=1,\dots, k$, let $a_r=\prod_{1\leq t\neq r \leq k} (\alpha_r-\alpha_t)$. Assume that the operator $M_z$ of multiplication by $z$ on the de Branges-Rovnyak space $H(B)$ is bounded. Then the Cauchy dual $M_z'$ of $M_z$ is subnormal if and only if the matrix 
	\[\displaystyle\sum_{r,t=1}^{k}\left( \frac{1}{a_r\overline{a_t}}\sum_{j=1}^{k}p_j(\alpha_r)\overline{p_j(\alpha_t)}\right)\left(1-\frac{1}{\alpha_r \overline{\alpha_t}}\right)^l\left(\left(\frac{1}{\alpha_r^{m+2}\overline{\alpha_t}^{n+2}}\right)\right)_{m,n\geq 0}\] 
	is formally positive semi-definite for every $l\geq 1$.
\end{theorem}
\begin{cor}\cite[Corollary 4.3]{cgr2022}\label{cor1}
	Assume the hypothesis of Theorem \ref{thm3}. If $\alpha_r\overline{\alpha_t}\notin [1,\infty)$ for every $1\leq r\neq t\leq k$, then $M_z'$ is subnormal if and only if 
	$\displaystyle \sum_{j=1}^k p_j(\alpha_r)\overline{p_j(\alpha_t)}=0, ~~~1\leq r\neq t\leq k.$
\end{cor}

As a first step towards the proof of the main theorem, we proceed with the following lemma. The lemma invokes the Riesz-Fejer theorem and derives a special case of equation (\ref{eq3}) by choosing $n=2, \zeta_1=1, \zeta_2=\xi$ and $c_1=c_2=1.$ Further in this case, we make some observations which turn out to be decisive in the proof of the main theorem. 

\begin{lemma}\label{l1}
	Let $\xi=\cos \theta +i \sin \theta\in \mathbb{T}~ (\xi\neq \pm1
    )$  be a point of the unit circle $\mathbb{T}.$
	\begin{enumerate}
		\item There exist a polynomial, say $f(z)$, of degree $4$ such that for all $z\in \mathbb{T},$ 
		\[|z-1|^2|z-\xi|^2+|z-1|^2+|z-\xi|^2=\displaystyle\frac{\overline{\xi}}{z^2}f(z).\]
		\item The polynomial $f(z)$ has exactly two roots lying outside the closed unit disc.
		\item Let $\alpha_1$ and $\alpha_2$ be those two roots of $f(z)$ lying outside the closed unit disc. Then $\alpha_1\alpha_2=b~\xi$ for some $b>0$.
		\item For all $z\in \mathbb{T}$, there exists a polynomial $q(z)$ of degree $2$ and a constant $d>0$ such that \begin{equation}\label{eq13}
			|z-1|^2|z-\xi|^2+|z-1|^2+|z-\xi|^2=d~ q(z)\overline{q(z)}.
		\end{equation}
		\item The sum $\alpha_1+\alpha_2=a(1+\xi)$ where $a=\frac{3b}{b+1}$.
		\item \label{part6} If $arg(\alpha_1)\neq \arg(\alpha_2)$ then $\alpha_2=\xi~\overline{\alpha_1}$ and thus $|\alpha_1|=|\alpha_2|$.
		\item \label{part7} The value of $b$ as obtained in (3) is a root of the polynomial \[g(x)=x^4-(6+2\cos\theta)x^3+(4+14\cos\theta)x^2-(6+2\cos\theta)x+1.\]
	\end{enumerate}
\end{lemma}
\ 
\begin{proof}
	\ 
	\begin{enumerate}
		\item After simplifying the expression $|z-1|^2|z-\xi|^2+|z-1|^2+|z-\xi|^2,$  
		%\begin{align*}
		%	&|z-1|^2|z-\xi|^2+|z-1|^2+|z-\xi|^2\\ &=\displaystyle\frac{\overline{\xi}z^4-3(1+\overline{\xi})z^3+(8+2\cos \theta)z^2-3(1+\xi)z+\xi}{z^2}\\&=\displaystyle\frac{\overline{\xi}f(z)}{z^2}
	%	\end{align*}
		we get, \begin{equation}\label{eq2}
			|z-1|^2|z-\xi|^2+|z-1|^2+|z-\xi|^2= \displaystyle\frac{\overline{\xi}f(z)}{z^2} 
		\end{equation}
		where $f(z)=z^4-3(1+\xi)z^3+(8+2\cos\theta)\xi z^2-3(\xi+\xi^2)z+\xi^2$.
		
		\item We know that $f(z)$ cannot have a root on the unit circle. Further, note that $f(z)$ is a monic polynomial with the property that if $\alpha$ is a root then $\frac{1}{\overline{\alpha}}$ is also a root. Therefore, exactly two of the roots of $f(z)$ lie outside the closed unit disc.
		
		\item The roots of $f(z)$ are $\alpha_1, \alpha_2, \frac{1}{\overline{\alpha_1}}$ and $\frac{1}{\overline{\alpha_2}}$. Therefore, we get,  \[ \alpha_1\alpha_2\displaystyle\frac{1}{\overline{\alpha_1}}\frac{1}{\overline{\alpha_2}}=\xi^2 \text{ and hence, } \displaystyle\frac{\alpha_1\alpha_2}{\overline{\alpha_1}\overline{\alpha_2}}=\frac{\xi}{\overline{\xi}}.\]
		This implies, 
		\begin{equation}\label{eq14}
			\alpha_1\alpha_2 \bar \xi=b \text{   for some } b\in \mathbb{R}. \text{ Hence }  \alpha_1\alpha_2 =b \xi. 
		\end{equation}
		In fact, we prove that $b>0$ and so that  $b=|\alpha_1||\alpha_2|$.\\
		Note that since left side of equation (\ref{eq2}) is positive, we have, $\frac{\bar \xi }{z^2}f(z)>0$, for all $z\in \mathbb{T}$.\\
		\noindent We define $q(z)$ as $q(z)=(z-\alpha_1)(z-\alpha_2)=z^2-(\alpha_1+\alpha_2)z+\alpha_1\alpha_2$. \\
	%	Let $\beta_1,\beta_2$ be the other two roots of $f(z)$. 
		Let, $\beta_1=\frac{1}{\overline{\alpha_1}}$ and $\beta_2=\frac{1}{\overline{\alpha_2}}.$
		 For $z\in \mathbb{T}$, we have
		\begin{align*}
			\displaystyle\frac{(z-\beta_1)(z-\beta_2)}{z^2}
		%	&=\displaystyle\frac{z^2-(\beta_1+\beta_2)z+\beta_1\beta_2}{z^2}\\
		%	&=1-\displaystyle\frac{\beta_1+\beta_2}{z}+\frac{\beta_1 \beta_2}{z^2}\\
			&=\beta_1\beta_2\left[\overline{z}^2-\left( \displaystyle\frac{1}{\beta_1}+\displaystyle\frac{1}{\beta_2}\right)\overline{z}+\displaystyle\frac{1}{\beta_1\beta_2} \right].
		\end{align*}
		But we have, $\frac{1}{\beta_1}=\overline{\alpha_1}$ and $\frac{1}{\beta_2}=\overline{\alpha_2}$. Thus, using equation (\ref{eq14}), we get,
		$\displaystyle\frac{(z-\beta_1)(z-\beta_2)}{z^2}=\displaystyle\frac{\xi}{b}\overline{q(z)}.$
		Hence, we get, 
		\begin{equation}\label{eq18}
			\displaystyle\frac{\bar \xi}{z^2}f(z)=\bar \xi \frac{\xi}{b}q(z)\overline{q(z)}=\frac{1}{b}\big|q(z)\big|^2,\text{ for all }z\in \mathbb{T}.
		\end{equation}
		This proves that $b>0$ and therefore $b=|\alpha_1||\alpha_2|$.
		
		\item We use equations (\ref{eq2}) and (\ref{eq18}) to get, 
		\[|z-1|^2|z-\xi|^2+|z-1|^2+|z-\xi|^2=d~ q(z)\overline{q(z)}\]
		where $d=\frac{1}{b}>0$.        
		\item We have $\alpha_1\alpha_2=b~\xi$ and  $\alpha_1+\alpha_2+\frac{1}{\overline{\alpha_1}}+\frac{1}{\overline{\alpha_2}}=3(1+\xi)$. 
		\\But, $\frac{1}{\overline{\alpha_1}}+\frac{1}{\overline{\alpha_2}}=\frac{\overline{\alpha_1}+\overline{\alpha_2}}{\overline{\alpha_1}~\overline{\alpha_2}}=\frac{\xi}{b}(\overline{\alpha_1}+\overline{\alpha_2})$. 
		\\Hence, $\alpha_1+\alpha_2+\frac{\xi}{b}(\overline{\alpha_1}+\overline{\alpha_2})=3(1+\xi)$.\\
		Let $\alpha_1+\alpha_2=r+is$, 
		then we get, $$b(r+is)+(\cos\theta+i\sin \theta)(r-is)=3b(1+\cos\theta+i\sin \theta).$$
		
		%     Thus, 
		%     $$br+r\cos\theta+s\sin\theta+(bs+r\sin\theta-s\cos\theta)i=3b(1+\cos\theta)+3b\sin\theta ~i.$$
		%     Equating real and imaginary parts, we get,
		% $$br+r\cos\theta+s\sin\theta=3b(1+\cos\theta) \text{ and } bs+r\sin\theta-s\cos\theta=3b\sin\theta.$$ 
		% which gives,
		% $$s\sin^2 \theta-s(b-\cos\theta)(b+\cos\theta)=3b(1+\cos\theta)\sin\theta-3b\sin\theta(b+\cos\theta).$$
		% Hence, we get, 
		% $$s(1-b^2)=3b\sin\theta+3b\sin\theta\cos\theta-3b^2\sin\theta-3b\sin\theta\cos\theta=3b\sin\theta(1-b)$$
		\noindent This gives, $s=a\sin\theta$  and $r=a(1+\cos\theta),$ where $a=\frac{3b}{b+1}.$
		% \\Substituting $s$ in one of the above equations, we get,
		% \\$r\sin\theta+a\sin\theta(b-\cos\theta)=3b\sin\theta$
		% \\$r=3b-ab +a\cos\theta=a+a\cos\theta=a(1+\cos\theta)$
		\\Thus, $\alpha_1+\alpha_2=a(1+\cos\theta+\sin\theta~i)=a(1+\xi)$.
		
		\item It can be verified that $f(\xi\overline{\alpha_1})=0$. We want to prove that $\alpha_2=\xi \overline{\alpha_1}$. 
		%If $\xi\overline{\alpha_1}=\frac{1}{\overline{\alpha_1}}$ then we get,  $|\alpha_1|=1$ which is not true. Further, if $\xi\overline{\alpha_1}=\frac{1}{\overline{\alpha_2}}$ then we get,  $|\alpha_1\alpha_2|=1$ which is not true. 
		Since, $|\alpha_1|,|\alpha_2|>1$, it is clear that $\xi\overline{\alpha_1}\neq \frac{1}{\overline{\alpha_1}}$ and $\xi\overline{\alpha_1}\neq \frac{1}{\overline{\alpha_2}}$.
		Hence, it remains to prove that $\xi\overline{\alpha_1}\neq \alpha_1$. 
		\\Suppose $\xi\overline{\alpha_1}= \alpha_1$. 
		Let $\arg(\alpha_i)=\theta_i$, for $i=1,2$. 
		We get, $\theta-\theta_1=\theta_1$ and thus, $\theta=2\theta_1$. 
		But,  $\alpha_1\alpha_2=b\xi$ gives $\theta_1+\theta_2=\theta$. Thus, $\theta_2=\theta_1$. But, this is not true. Hence, $\xi \overline{\alpha_1}\neq \alpha_1$.\\
		Thus, we must have $\alpha_2=\xi \overline{\alpha_1}$ and hence $|\alpha_1|=|\alpha_2|$.

		\item The relation between the roots and the coefficients of $f(z)$ gives,
		\[\alpha_1\alpha_2+\alpha_1\frac{1}{\overline{\alpha_1}}+\alpha_1\frac{1}{\overline{\alpha_2}}+\alpha_2\frac{1}{\overline{\alpha_1}}+\alpha_2\frac{1}{\overline{\alpha_2}}+\frac{1}{\overline{\alpha_1}~\overline{\alpha_2}}=(8+2\cos \theta)\xi.\]
		Using the above relation, it can be seen that  
		%		$1+|a(1+\xi)|^2=8b-b^2+2b\cos \theta$   and hence,
	%	$1+\displaystyle\frac{9b^2}{(b+1)^2}\langle1+\xi,1+\xi\rangle=8b-b^2+2b\cos\theta.$\\
	%	Hence
	%	\\$b^4-2(3+\cos\theta)b^3+2(2+7\cos\theta)b^2-2(3+\cos\theta)b+1=0.$
		\noindent $b$ satisfies the polynomial 
		$$g(x)=x^4-(6+2\cos\theta)x^3+(4+14\cos\theta)x^2-(6+2\cos\theta)x+1.$$
It is easy to observe that for a non-zero real number $\beta$, if $g(\beta)=0$ then $g(\frac{1}{\beta})=0.$
	\end{enumerate}

\end{proof}

The proof of the main theorem demands different techniques based on the location of the point $\xi=\cos \theta + i\sin \theta$ on the unit circle. We therefore divide the proof into three parts. In fact, the nature of the roots $\alpha_1, \alpha_2$ appearing in Lemma \ref{l1} when $\cos \theta =0.6$ prompts the authors to make these cases, as can be seen in the following lemma.   
\begin{lemma}\label{l3}
	Let $\xi=\cos \theta + i\sin \theta \in \mathbb{T}$ $(\xi\neq \pm 1)$. Let $\alpha_1,\alpha_2$ be two complex numbers outside closed unit disc, obtained in part (3) of Lemma \ref{l1}.
	\begin{enumerate}
		\item If $\cos\theta<0.6$ then $\arg(\alpha_1)\neq \arg(\alpha_2)$.
		\item If $\cos\theta>0.6$ then $\arg(\alpha_1)=\arg(\alpha_2)$ and $|\alpha_1|\neq |\alpha_2|.$
		\item If $\cos\theta=0.6$ then $\alpha_1=\alpha_2=2+i$.
	\end{enumerate}
\end{lemma}
\begin{proof}
	\ 
	\begin{enumerate}
		\item We have $\cos\theta<0.6$. Suppose 		$\arg(\alpha_1)=\arg(\alpha_2)$. \\
		Now as $\alpha_1+\alpha_2=a(1+\xi)$ we have $\alpha_1=s(1+\xi)$ and $\alpha_2=t(1+\xi)$ for some $s,t\in \mathbb{R}$.
		This gives $\alpha_1+\alpha_2=(s+t)(1+\xi)=a(1+\xi)$ and $\alpha_1\alpha_2=st(1+\xi)^2=b\xi$.\\
		Thus we get, $s+t=a$ and $st=\displaystyle\frac{b\xi}{(1+\xi)^2}=\frac{b}{|1+\xi|^2}=\frac{b}{2(1+\cos\theta)}$.\\
		
		Consider a polynomial $h(x)=(x-s)(x-t)$. 
		\\Then $h(x)=x^2-ax+\displaystyle\frac{b}{2(1+\cos\theta)}=x^2-\displaystyle\frac{3b}{b+1}x+\displaystyle\frac{b}{2(1+\cos\theta)}$.\\
		Let $D$ be the discriminant of polynomial $h(x)$ then,
		\[D=\displaystyle\frac{9b^2}{(b+1)^2}-\frac{2b}{1+\cos\theta}.\]
		Since $\cos\theta<0.6$, we have $\displaystyle\frac{-1}{1+\cos\theta}<\frac{-1}{1.6}$.\\
	%	Thus,
	%	\begin{align*}
	%		D&=\displaystyle\frac{9b^2}{(b+1)^2}-\frac{2b}{1+\cos\theta}\\
	%		&<\displaystyle\frac{9b^2}{(b+1)^2}-\frac{2b}{1.6}\\
	%		&=\displaystyle\frac{9b^2}{(b+1)^2}-\frac{5b}{4}\\
	%		&= \displaystyle\frac{36b^2-5b(b^2+2b+1)}{4(b+1)^2}\\
	%		&=\displaystyle\frac{-b(5b^2-26b+5)}{4(b+1)^2}\\
	%		&=\displaystyle\frac{-b(5b-1)(b-5)}{4(b+1)^2}
	%	\end{align*}
		It can be seen that  $D<\displaystyle\frac{-b(5b-1)(b-5)}{4(b+1)^2}$.\\
		Note that $b>1$ and hence we have $5b-1>0$.\\
		We now prove that $b>5$.\\
		Recall from part (\ref{part7}) of Lemma \ref{l1} that $b$ satisfies the polynomial \[g(x)=x^4-(6+2\cos\theta)x^3+(4+14\cos\theta)x^2-(6+2\cos\theta)x+1.\]
		We further note that 
		\begin{equation}\label{eq3.18.2}
			\begin{aligned}
				g(0)=1,\\
				g(1)=-6+10 \cos\theta,\\
				g(3)=-62+66\cos\theta,\\
				g(4)=-87+88\cos\theta,\\
				g(5)=-54+90\cos\theta,\\
				g(6)=109-60\cos\theta.
			\end{aligned}
		\end{equation}
		Clearly $g(0)>0$. Also, since, $\cos\theta<0.6$, using (\ref{eq3.18.2}), we get, 
		\[g(1)<0, g(5)<0 \text{ and }g(6)>0.\] 
		Thus there exists $\beta_0\in (5,6)$ such that $g(\beta_0)=0$.  %implying $g\left(\frac{1}{\beta_0}\right)=0.$ 
		We now prove that $g$ has no root in the interval $(1,5)$. Suppose  $g$ has a root in the interval $(1,5)$. Let $\beta$ be the smallest such a root of $g.$
		\begin{enumerate}
			\item [Case (i)] Let $g(x)<0$, $\forall x\in (\beta,5)$.\\
			Then $\beta$ is a local maxima of $g$. Hence, $g'(\beta)=0$.\\
			But, $g$ is a polynomial with $g(\beta)=0$ and $g'(\beta)=0$. Thus, $\beta$ is a repeated root of $g$. Now by the property of $g$, it has one more repeated root (that is , $\frac{1}{\beta}$). Thus, $g$ has more than $4$ roots which is a contradiction.
			\item [Case (ii)] Suppose there exists $x_0\in (\beta,5)$ such that $g(x_0)>0$.\\
			But, $g(5)<0$ and hence there exists $\gamma\in (x_0,5)$ such that $g(\gamma)=0$. Thus, we get, $g(\frac{1}{\gamma})=0$. This implies that $g$ has more than $4$ roots, which is a contradiction.
		\end{enumerate}
		Hence, $g$ has no root in the interval $(1,5)$.
		Therefore, $b>5$.\\
		Thus, if $\cos\theta<0.6$ then $D<0$ contradicting the assumption that  $s,t$ are the real roots of the polynomial $h(x)$. Hence,	\\$\arg(\alpha_1)\neq\arg(\alpha_2).$
		\item Let $\cos\theta>0.6.$ 
		We need to check whether there exist $s,t\in \mathbb{R}$ such that $\alpha_1=s(1+\xi)$ and $\alpha_2=t(1+\xi)$. Indeed, it is enough to show that there exist $s,t\in \mathbb{R}$ such that 		 
		$\alpha_1+\alpha_2=(s+t)(1+\xi)=a(1+\xi)$ and $\alpha_1\alpha_2=st(1+\xi)^2=b\xi$.\\ Following the steps in (1), we have 
	%	Thus we get $s+t=a$ and $st=\displaystyle\frac{b\xi}{(1+\xi)^2}=\frac{b}{|1+\xi|^2}=\frac{b}{2(1+\cos\theta)}$.\\
		%Consider a polynomial $h(x)=(x-s)(x-t)$. 
	%	\\Then $h(x)=x^2-ax+\displaystyle\frac{b}{2(1+\cos\theta)}=x^2-\displaystyle\frac{3b}{b+1}x+\displaystyle\frac{b}{2(1+\cos\theta)}$.\\
		the discriminant $D$ of polynomial $h(x)$ as
		\[D=\displaystyle\frac{9b^2}{(b+1)^2}-\frac{2b}{1+\cos\theta}.\]
		Here, we have, $\cos\theta>0.6$, which gives $\displaystyle\frac{-1}{1+\cos\theta}>\frac{-1}{1.6}$.\\
	%	Thus,
	%	\begin{align*}
	%		D&=\displaystyle\frac{9b^2}{(b+1)^2}-\frac{2b}{1+\cos\theta}\\
	%		&>\displaystyle\frac{9b^2}{(b+1)^2}-\frac{2b}{1.6}=\displaystyle\frac{9b^2}{(b+1)^2}-\frac{5b}{4}\\
	%		&= \displaystyle\frac{36b^2-5b(b^2+2b+1)}{4(b+1)^2}=\displaystyle\frac{-b(5b^2-26b+5)}{4(b+1)^2}\\
	%		&=\displaystyle\frac{-b(5b-1)(b-5)}{4(b+1)^2}
	%	\end{align*}
		Thus, $D>\displaystyle\frac{-b(5b-1)(b-5)}{4(b+1)^2}.$\\
		Note that $b>1$ and hence we have $5b-1>0$.\\
		Now, for $\cos\theta>0.6$, from equation (\ref{eq3.18.2}), we get
		\[g(1)>0, g(5)>0, g(6)>0.\]
%		Note that $g(4)<0$ if $\cos\theta <\displaystyle\frac{87}{88}$ and
%		$g(3)<0$ if $\cos\theta<\displaystyle\frac{62}{66}$.\\
%		Thus, if $0.6<\cos\theta<\frac{62}{66}$ then $g(3)<0$ which implies that there exist $\beta_0\in (1,3)$ and $\gamma_0\in (3,5)$ such that $g(\beta_0)=0=g(\gamma_0)$. Thus we have  $g(\frac{1}{\beta_0})=0=g(\frac{1}{\gamma_0})$.
%		This implies that all the roots of $g$ are less than $5$.\\
%		Now, if $\frac{62}{66}<\cos\theta<\frac{87}{88}$ then $g(4)<0$. A similar argument as above shows that all the roots of $g$ are less than $5$. \\
%		Finally, if $\frac{87}{88}<\cos\theta<1$  then note that $g(3)>0$ as well as $g(4)>0$.\\
		It can be seen that $g(2+\sqrt{3})=-42-24\sqrt{3}+42\cos\theta+24\sqrt{3}\cos\theta$.\\
		Thus, $g(2+\sqrt{3})<0$ since $\cos\theta<1$. 
		This assures that $g$ has two roots in the interval $(1,5)$ and hence another two roots in the interval $(0,1)$.
		Thus, all the roots of $g$ are less than $5$.
		Hence, we have $b<5$ implying $D>0$.
		Hence, there exist distinct real numbers $s, t$ such that $\alpha_1=s(1+\xi)$ and $\alpha_2=t(1+\xi)$.
		Therefore, $\arg(\alpha_1)=\arg(\alpha_2)$.\\
		Further, note that since $s,t$ are distinct real roots, we must have $|\alpha_1|\neq |\alpha_2|$.
		\item It can be verified that if $\cos\theta=0.6$ then $2+i$ is a repeated root of $f(z)$. Therefore, $\alpha_1=\alpha_2=2+i$.
	\end{enumerate}
\end{proof}

The proof of the main theorem relies on the technique discussed in \cite{cgr2022}. The computational details in this case have been discussed below.  \\
The polynomial $q(z)$ obtained in the equation (\ref{eq13})  allows us to compute $O_\mu(z)$ as follows:
$O_\mu(z)=\displaystyle\frac{p(z)}{q(z)}$ where $p(z)=\displaystyle\frac{e^{i\theta_0}}{\sqrt{d}}(z-1)(z-\xi)$ and \\ $q(z)=(z-\alpha_1)(z-\alpha_2)$.
\\ Now $O_\mu(0)=\displaystyle\frac{p(0)}{q(0)}=\frac{e^{i\theta_0}(-1)(-\xi)}{\sqrt{d}\alpha_1\alpha_2}=\frac{e^{i\theta_0}}{\sqrt{d}~b}$.
\\Hence, we choose $\theta_0=0$ to have $O_\mu(0)>0$. \\ Thus, we have, 
\begin{equation}\label{eq21}
	p(z)=\displaystyle\frac{1}{\sqrt{d}}(z-1)(z-\xi) \text{ and }
	O_\mu(z)=\displaystyle\frac{1}{\sqrt{d}}\frac{(z-1)(z-\xi)}{q(z)}.
\end{equation}

A close look at \cite[Theorem 6.4]{cgr2022} suggests that we need to first compute the matrix $\left(\left[\langle f_i,f_j\rangle_\mu\right]_{1\leq i,j\leq 2}\right)^{-1}$ (say $B$) in the equation (\ref{eq7}). We thus proceed for the computations of the entries of this matrix through the following lemma.
\begin{lemma}\label{l2}
	Let matrix $C=\left[c_{ij}\right]_{1\leq i,j\leq 2}$, where $c_{ij}=\langle f_i,f_j\rangle_\mu$, $1\leq i,j\leq 2$ and $f_1,f_2$ as defined in (\ref{eq11}). Then 
	\begin{enumerate}
		\item $\Delta=\det (C)=\left[(a-3)\displaystyle\frac{1+\cos\theta}{1-\cos\theta}+\frac{3}{a(1-\cos\theta)}\right]$.
		\item Matrix $B$ is given by $B=\begin{pmatrix}
			a-1 & -\displaystyle\frac{c_{12}}{\Delta}\\ -\displaystyle\frac{\overline{c_{12}}}{\Delta} & a-1
		\end{pmatrix},$
		where $c_{12}$ is given by $c_{12}=\displaystyle\frac{d q(1)^2}{(\xi-1)^3}$.
	\end{enumerate}
\end{lemma}
\begin{proof}

	\begin{enumerate}
		\item We have, $O_\mu(z)=\displaystyle\frac{1}{\sqrt{d}}\frac{(z-1)(z-\xi)}{q(z)}$.
		\\Thus, $O_\mu'(1)=\displaystyle\frac{1}{\sqrt{d}}\frac{1-\xi}{q(1)}$ and $O_\mu'(\xi)=\displaystyle\frac{1}{\sqrt{d}}\frac{\xi-1}{q(\xi)}.$
		\\
		
		\noindent We observe that $|O_\mu'(1)|^2=1=|O_\mu'(\xi)|^2$ and $q(\xi)=\xi^2~\overline{q(1)}$.\\
		
		\noindent This gives the relation : $O_\mu'(\xi)=\overline{\xi}~\overline{O_\mu'(1)}.$
		\\We use equation (\ref{eq11}) to define $f_j(z)=\displaystyle\frac{O_\mu(z)}{O_\mu'(\xi_j)(z-\xi_j)}$.
		\\Thus, \begin{equation}\label{eq24}
			f_1(z)=\displaystyle\frac{q(1)}{1-\xi}\frac{z-\xi}{q(z)} \text{ and } f_2(z)=\displaystyle\frac{z-1}{q(z)}\frac{q(\xi)}{\xi-1}.
		\end{equation}
		% \\$f_1'(z)=\displaystyle\frac{q(1)}{1-\xi}\left[  \frac{q(z)-(z-\xi)q'(z)}{q(z)^2}\right]$.
		
		\noindent This gives, $f_1'(1)=\displaystyle\frac{1}{1-\xi}-\frac{q'(1)}{q(1)}$.\\
		
		\noindent Hence, by equation (\ref{eq8}) we get, $c_{11}=f_1'(1)=\displaystyle\frac{1}{1-\xi}-\frac{q'(1)}{q(1)}.$ \\
		But, $q(z)=z^2-a(1+\xi)z+b\xi$, hence, we get, 
	%	\begin{align*}
	%		c_{11}&=\displaystyle\frac{1}{1-\xi}+\frac{a-2+a\xi}{1-a+(b-a)\xi}\\
			% &=\displaystyle\frac{1}{1-\xi}+\frac{[(a-2)+a\xi][(1-a)+(b-a)\overline{\xi}]}{[(1-a)+(b-a)\xi][(1-a)+(b-a)\overline{\xi}]}\\
			% &=\displaystyle\frac{1}{1-\xi}+\frac{(a-2)(1-a)+a(b-a)+\xi a(1-a)+(a-2)(b-a)\overline{\xi}}{x^2+y^2+2xy\cos \theta}\\
	%		&=\Re\left( \displaystyle\frac{1}{1-\xi}\right)+\frac{(a-2)(1-a)+a(b-a)+a(1-a)\cos \theta +(a-2)(b-a)\cos \theta}{x^2+y^2+2xy\cos \theta},
	%	\end{align*}
	%	where $x=1-a$ and $y=b-a$.
	%	Thus we get,
	%	\[x^2+y^2+2xy\cos\theta=\displaystyle\frac{g(b)+2(1-\cos\theta)b(b^2+2b+1)}{(b+1)^2}=2b(1-\cos\theta)\]
		 \[c_{11}=\displaystyle\frac 12+\frac{(a-2)(1-a)+a(b-a)+a(1-a)\cos\theta+(a-2)(b-a)\cos\theta}{2b(1-\cos\theta)}.\]
		
		\noindent It can be seen that the numerator of $c_{11}$ has $(a-1)$ as a factor, this along with the observation that $b(3-a)=a$ gives, 
		\begin{equation}\label{eq19}
			c_{11}=(a-1)\left[(a-3)\displaystyle\frac{1+\cos\theta}{1-\cos\theta}+\frac{3}{a(1-\cos\theta)}\right].
		\end{equation}
		We once again use equations (\ref{eq8}) and (\ref{eq11}) to get, 
		$$c_{22}=\xi f_2'(\xi)
		\text{ where } f_2(z)=\displaystyle\frac{(z-1)}{\sqrt{d}q(z)O_\mu'(\xi)}.$$
	%	We further have, $\displaystyle\frac{O_\mu(z)}{z-\xi}=\frac{1}{\sqrt{d}}\frac{z-1}{q(z)}.$
		It is easy to see that  $f_2'(z)=\displaystyle\frac{1}{\sqrt{d}O_\mu'(\xi)}\left[ \frac{q(z)-(z-1)q'(z)}{q(z)^2}\right].$
		% \\Thus, $f_2'(\xi)=\displaystyle\frac{1}{\sqrt{d}O_\mu'(\xi)}\left[ \frac{q(\xi)-(z-1)q'(\xi)}{q(\xi)^2}\right]=\frac{\sqrt{d}q(\xi)}{(\xi-1)\sqrt{d}}\left[ \frac{1}{q(\xi)}-\frac{(\xi-1)q'(\xi)}{q(\xi)^2}\right].$
		\\Hence, $f_2'(\xi)=\displaystyle\frac{1}{\xi-1}-\frac{q'(\xi)}{q(\xi)}.$
		Thus, $c_{22}=\xi f_2'(\xi)=\displaystyle\frac{\xi}{\xi-1}-\frac{\xi q'(\xi)}{q(\xi)}.$
		%\\Note that $\displaystyle\frac{\xi}{\xi-1}=\frac{1}{1-\overline{\xi}}$, $\overline{\xi}q'(\xi)=\overline{q'(1)}$.
		It can be seen that
		$$c_{22}=\displaystyle\frac{1}{1-\overline{\xi}}-\frac{\xi q'(\xi)}{\xi^2\overline{q(1)}}=\frac{1}{1-\overline{\xi}}-\frac{\overline{q'(1)}}{\overline{q(1)}}=\overline{c_{11}}.$$
		\\Now we use equation (\ref{eq9}) to obtain the following:
		$$c_{12}=\displaystyle\frac{1}{O_\mu'(1)O_\mu'(\xi)(1-\overline{\xi})}
		% =\frac{1}{\xi O_\mu'(1)^2(1-\overline{\xi})}=\frac{1}{(\xi-1)O_\mu'(1)^2}=\frac{1}{\xi-1}\frac{d~q(1)^2}{(1-\xi)^2}
		=\frac{d q(1)^2}{(\xi-1)^3}.$$
		Thus, after simplification we get,  $|c_{12}|^2
		% =\displaystyle\frac{d^2|q(1)|^4}{|\xi-1|^6}=\frac{d^2|1-\xi|^4}{d^2|1-\xi|^6}=\frac{1}{|1-\xi|^2}
		=\displaystyle\frac{1}{2(1-\cos\theta)}.$\\
				\noindent Observe that $c_{21}=\overline{c_{12}}.$\\
		We express $\Delta=\det(C)=c_{11}c_{22}-|c_{12}|^2$ as follows: 
        \begin{align*}
			\Delta&=c_{11}c_{22}-|c_{12}|^2\\
			&=\displaystyle\frac{a(a-3)(1+\cos \theta )+3}{a(1-\cos\theta)}\Bigg[ \frac{(a-1)^2(a(a-3)(1+\cos\theta)+3)}{a(1-\cos\theta)}\\
			&\hspace{4.5cm}-\frac{a}{2(a(a-3)(1+\cos\theta)+3)}\Bigg].
		\end{align*}
		\noindent We now prove that $$\left[ \frac{(a-1)^2(a(a-3)(1+\cos\theta)+3)}{a(1-\cos\theta)}
		-\frac{a}{2(a(a-3)(1+\cos\theta)+3)}\right]=1.$$
		\noindent It is enough to prove that 
		\[2(a-1)^2(a(a-3)(1+\cos\theta)+3)^2-a^2(1-\cos\theta)-2a(1-\cos\theta)(a(a-3)(1+\cos\theta)+3)=0.\]
		But, we know that $a-1=\frac{2b-1}{b+1}$, and thus the above equation reduces to  
		$9(3b^2\cos\theta+5b^2-6b\cos\theta-2b+2)g(b)=0.$\\
		This follows from the fact that $b$ is a root of the polynomial $g(z).$
		\\Hence, 
		\begin{equation}\label{eq20}
			\Delta=\det (C)=\left[(a-3)\displaystyle\frac{1+\cos\theta}{1-\cos\theta}+\frac{3}{a(1-\cos\theta)}\right].
		\end{equation}
		\item The equations (\ref{eq19}) and (\ref{eq20}) imply that $b_{11}=c_{11}/\Delta = a-1$.\\
		Also, as $c_{22}=\overline{c_{11}}$, we get, $b_{22}=c_{22}/\Delta=a-1$.\\
		Thus, we get, \begin{equation}\label{eq25}
			B=\begin{pmatrix}
				a-1 & -\displaystyle\frac{c_{12}}{\Delta}\\ -\displaystyle\frac{\overline{c_{12}}}{\Delta} & a-1
			\end{pmatrix},
		\end{equation}
		where $c_{12}=\displaystyle\frac{d q(1)^2}{(\xi-1)^3}$.
	\end{enumerate}
\end{proof}

We now derive the following observations about $b$, $a$ and $\Delta$ as functions of $\theta$ when $\cos\theta <0.6.$ These observations are used in the proof of the main theorem.

\begin{lemma}\label{l6}
Let $\xi=\cos \theta + i\sin \theta \in \mathbb{T}$ $(\xi\neq \pm 1)$. Let $\alpha_1,\alpha_2$ be two complex numbers outside closed unit disc, as obtained in part (3) of Lemma \ref{l1}. Let $b=|\alpha_1\alpha_2|$, $a=\frac{3b}{b+1}$ and $\Delta$ as obtained in part (1) of Lemma \ref{l2}. Let $t_0=\cos^{-1}(0.6)$. Then, for $\theta\in [t_0,\pi/2],$
\begin{enumerate}
    \item $b$ is an increasing function of $\theta.$
    \item $a$ is an increasing function of $\theta.$
    \item $0.7<\Delta<1.10056.$
\end{enumerate}
\end{lemma}
\begin{proof}
\ 
    \begin{enumerate}
        \item We note from part (7) of Lemma \ref{l1} that $b$ is a root of \\$g(x)=x^4-Ax^3+Bx^2-Ax+1,$ where $A=6+2\cos\theta, B=4+14\cos\theta.$\\
    Thus, $\frac{g(x)}{x^2}=h(y),$ where $h(y)=y^2-Ay+B-2,$ where $y=x+\frac{1}{x}.$\\
    $h(y)=0\implies y=\frac{A\pm \sqrt{A^2-4B+8}}{2}.$\\
    It can be seen that $\frac{A+ \sqrt{A^2-4B+8}}{2}>2.$ Let $y_1=\frac{A+ \sqrt{A^2-4B+8}}{2}>2.$\\
    Let $c(\theta)=\cos\theta$ be a function defined on $[t_0,\pi/2].$\\
    Then we get,\\
    $y_1(c)=3+c+\sqrt{c^2-8c+7}$ and hence $\frac{dy_1}{dc}=1+\frac{c-4}{\sqrt{c^2-8c+7}}.$\\
    Now, $|c-4|>\sqrt{c^2-8c+7}$ and $c<4$ gives, $\frac{dy_1}{dc}<0$ and hence $y_1$ is a decreasing function of $c$. Further, $c$ is  a decreasing function of $\theta$ and hence, $y_1$ is an increasing function of $\theta.$ This proves that $\frac{y_1(\theta)+\sqrt{y_1(\theta)^2-4}}{2}$ is an increasing function of $\theta.$ \\
    We know that $b$ is a root of $y=x+\frac{1}{x}.$ \\Thus we have $b(\theta)=\frac{y_1(\theta)\pm\sqrt{y_1(\theta)^2-4}}{2}.$ \\
    As $b>5,$ we must have $b(\theta)=\frac{y_1(\theta)+\sqrt{y_1(\theta)^2-4}}{2}.$ Hence, $b$ is an increasing function of $\theta.$
    \item We consider $a$ as a function of $b$ given by  $a(b)=\frac{3b}{b+1}$. It can be seen that $a'(b)>0$ and hence $a$ is an increasing function of $b$. Thus, $a$ is an increasing function of $\theta$ on $[t_0,\pi/2].$
    \item We know that for $\theta=t_0,$ $b=5$ (See Part (3) of Lemma \ref{l3}) and for $\theta=\frac{\pi}{2},$ $b=5.46269$ (Refer \cite[Appendix A.2]{mkvms2025}). Since, $a$ is an increasing function of $b$, we get, $a\in [2.5,2.5357]$ whenever \\$\theta\in [t_0,\pi/2].$\\
    We have $\Delta=(a-3)\displaystyle\frac{1+\cos\theta}{1-\cos\theta}+\frac{3}{a(1-\cos\theta)}.$\\
    For fixed $a\in [2.5,2.5357]$, we have, $\Delta'(\theta)=\frac{-(2a^2-6a+3)\sin\theta}{a(1-\cos\theta)^2}<0.$\\
    Thus, for fixed $a$, $\Delta$ is a decreasing function of $\theta.$\\
    For fixed $\theta\in [t_0,\pi/2]$, we have, $\Delta'(a)=\frac{a^2(1+\cos\theta)-3}{a^2(1-\cos\theta)}>0.$\\
    Thus, for fixed $\theta$, $\Delta$ is an increasing function of $a.$
    Now, for any $\theta\in [t_0,\pi/2]$ and $a\in [2.5,2.5357]$, we get,\\
    $\Delta(\theta,a)>\Delta(\theta,2.5)>\Delta(\pi/2,2.5)=0.7.$ and\\
    $\Delta(\theta,a)<\Delta(\theta,2.5357)<\Delta(t_0,2.5357)=1.10056.$\\
    Thus, we get, $0.7<\Delta<1.10056.$
    \end{enumerate}
\end{proof}

The following theorem is a generalization of \cite[Theorem 8]{mkvms2025}, developed in the case of $\{1,i\}$. We use proof of Theorem \ref{thm2} to find $B=(b_1,b_2)$ such that $D(\mu)$ coincides with $H(B)$ with equality of norms.
\begin{theorem}\label{thm5}
	Let $\xi=\cos \theta + i\sin \theta \in \mathbb{T}$ $(\xi\neq \pm 1)$. Let $\mu$ be a finite, positive, Borel measure on the unit circle $\mathbb{T}$ of the form $\mu= \delta_1+\delta_{\xi} $. Then the Dirichlet-type space $D(\mu)$ coincides with the de Branges-Rovnyak space $H(B)$ with equality of norms where $B=(b_1,b_2)$ and $b_j=\displaystyle\frac{p_j}{q}$, $j=1,2$ where 
	$$p_1(z)=p_{11}z+p_{12}z^2,~p_2(z)=p_{22}z^2 \text{ and } q(z)=(z-\alpha_1)(z-\alpha_2)$$
	with $\alpha_1,\alpha_2$ as obtained in the Lemma \ref{l1} and $p_{11}, p_{22}$ and $p_{12}$ are constants.
\end{theorem}
\begin{proof}
	Let $\xi_1=1$, $\xi_2=\xi$, $B=(b_{ij})_{1\leq i,j\leq 2}$ be the matrix obtained in Lemma \ref{l2}. Also, let $O_\mu(z)$ be defined by equation (\ref{eq21}).\\ 
	The operator $M_z$ on $D(\mu)$ is analytic, norm increasing (see \cite[Theorem 3.6]{richter1991}), and the reproducing kernel for $D(\mu)$ is normalized. Therefore, by \cite[Lemma 3.4]{cgr2022}, there exists a positive semi-definite kernel $\eta:\mathbb{D}\times \mathbb{D}\to \mathbb{C}$ such that 
	\[\eta(z,0)=0,~~~k(z,w)=\displaystyle\frac{1-\eta(z,w)}{1-z\overline{w}},~~z,w\in \mathbb{D}.\]
	We rewrite the expression for $K(z,w)$ using equations (\ref{eq4}), (\ref{eq5}) and \\Theorem \ref{thm1} as 
	\begin{align*}
		k(z,w)&=\displaystyle\frac{O_\mu(z)\overline{O_\mu(\lambda)}}{1-z\overline{w}}+\sum_{j=1}^{2}f_j(z)g_j(w)\\
		&=\displaystyle \frac{1}{1-z\overline{w}}\frac{p(z)\overline{p(w)}}{q(z)\overline{q(w)}}+\sum_{i,j=1}^{2}\overline{b_{ji}}f_j(z)\overline{f_i(w)}.
	\end{align*}	
	\noindent Equating both the above expressions for the reproducing kernel $k(z,w)$, we get, \\
	$q(z)\eta(z,w)\overline{q(w)}$ \\ $=q(z)\overline{q(w)}-p(z)\overline{p(w)}\left(1+(1-z\overline{w})\displaystyle\sum_{i,j=1}^{2}\frac{\overline{b_{ji}}}{O_\mu'(\xi_j)\overline{O_\mu'(\xi_i)}}\frac{1}{(z-\xi_j)(\overline{w}-\overline{\xi_i})}\right).$
	
	\noindent Since, the right hand side is a polynomial in $z$ and $\overline{w}$, there exists a matrix $\hat A=(a_{ij})_{0\leq i,j\leq 2}$ such that the above expression is equal to $\displaystyle\sum_{i,j=0}^{2}a_{ij}z^i\overline{w}^j$, $z,w\in \mathbb{D}$. The matrix $\hat A$ is positive semi-definite.\\
	As this expression at $w=0$ is $0$ for each $z\in \mathbb{D}$, we obtain a matrix $A$ from $\hat A$ by deleting first row and first column of matrix $\hat A$. Therefore, matrix $A$ is also positive semi-definite.\\
	We present here the computations to obtain matrix $A$.
	We begin with simplification of the right hand side.\\
	$q(z)\overline{q(w)}-p(z)\overline{p(w)}-p(z)\overline{p(w)}(1-z\overline{w})\displaystyle\sum_{i,j=1}^2\frac{\overline{b_{ji}}}{O_\mu'(\zeta_j)\overline{O_\mu(\zeta_i)}}\frac{1}{(z-\zeta_j)(\overline{w}-\overline{\zeta_i)}}$
	\begin{align*}
		&=(z^2-a(1+\xi)z+b\xi)(\overline{w}^2-a(1+\overline{\xi})\overline{w}+b\overline{\xi})\\&-\displaystyle\frac{1}{d}(z^2-a(1+\xi)z+b\xi)(\overline{w}^2-a(1+\overline{\xi})\overline{w}+b\overline{\xi})
		\\&-\frac{1}{d}(1-z\overline{w})\Bigg[\overline{b_{11}}(z\overline{w}-z\overline{\xi}-\xi\overline{w}+1)+\overline{b_{22}}(z\overline{w}-z-\overline{w}+1)
		\\&\hspace{2.3cm}+\frac{\overline{b_{12}}(z\overline{w}-z-\xi\overline{w}+\xi)}{\xi O_\mu'(1)^2}+\frac{\overline{b_{21}}(z\overline{w}-z\overline{\xi}-\overline{w}+\overline{\xi})}{\xi O_\mu'(\xi)^2} \Bigg].
	\end{align*}
	We get, $a_{11}=\text{coefficient of }z\overline{w}=$
	$(a^2-b)(|1+\xi|^2)+\displaystyle\frac{2b}{\Delta}.$
	
	$a_{22}=\text{coefficient of }z^2\overline{w}^2$  
	$=\displaystyle1-3b+2ab-\frac{b}{\Delta}.$
	\\$a_{12}=\text{coefficient of }z\overline{w}^2=-b(1+\xi).$  \\
	Note that $a_{11}$ and $det(A)=a_{11}a_{22}-|a_{12}|^2$ are non-negative real numbers. Thus,  $a_{22}$ is also a non-negative real number.\\
	%\[a_{12}
	%(1+\xi)(b-a)-(a-1)b(1+\xi)+0
	% =(1+\xi)[b-a-b(a-1)]\\
	% =(1+\xi)(2b-a-ab)=(1+\xi)(2b-3b)
%	=-b(1+\xi).
	As matrix $A=(a_{ij})_{1\leq i,j\leq 2}$ is positive semi-definite, by applying Cholesky's decomposition we get an upper triangular matrix $P$ such that $A=P^*P$, where 
	\[P=\begin{pmatrix}
		\sqrt{a_{11}} &\displaystyle\frac{a_{12}}{\sqrt{a_{11}}}\\
		0& \displaystyle \sqrt{\frac{a_{11}a_{22}-|a_{12}|^2}{a_{11}}}
	\end{pmatrix}.\]
	Once again by using Theorem \ref{thm2}, we get,
	\begin{align}
		\displaystyle p_1(z)&=\langle P (z,z^2)^T,(1,0)\rangle= \sqrt{a_{11}}  ~z+\frac{a_{12}}{\sqrt{a_{11}}} z^2, \label{eq22}\\
		\displaystyle p_2(z)&=\langle P (z,z^2)^T,(0,1)\rangle=p_{22}z^2=\sqrt{\frac{a_{11}a_{22}-|a_{12}|^2}{a_{11}}}~z^2\label{eq23}.
	\end{align}
	
	Therefore, we conclude that the Dirichlet-type space $D(\mu)$ coincides with de Branges-Rovnyak space $H(B)$ with equality of norms where $B=(b_1,b_2)$ and $b_j=\displaystyle\frac{p_j}{q}$, $j=1,2.$
\end{proof}
These computations lead to the proof of the main theorem of this paper. We consider $M_z$ on the above $D(\mu)$ and claim that its Cauchy dual is not subnormal.\\

{\bf Proof of the Main Theorem :}
\begin{proof}
{\bf Case 1:} Let $\cos\theta<0.6.$\\
	Note that if $\Im(\alpha_1\overline{\alpha_2})=0$ then $ \frac{\alpha_1}{\alpha_2}=\frac{\overline{\alpha_1}}{\overline{\alpha_2}}.$ This implies that $ \frac{\alpha_1}{\alpha_2}$ is a real number, which is not true as $\arg(\alpha_1)\neq \arg(\alpha_2)$. Thus, $\Im(\alpha_1\overline{\alpha_2})\neq 0.$\\
	Also, by Lemma \ref{l3}, we know that $\alpha_2=\xi \overline{\alpha_1}$.\\
    Thus, we have, $\alpha_1\overline{\alpha_2}\in \mathbb{C}\setminus[1,\infty)$ and $\alpha_2\overline{\alpha_1}\in \mathbb{C}\setminus[1,\infty)$. \\
	We now consider the expression:
	$p_1(\alpha_1)\overline{p_1(\alpha_2)}+p_2(\alpha_1)\overline{p_2(\alpha_2)}.$ \\
	Substituting, $\alpha_1=\xi \overline{\alpha_2}$ in the above expression, we get \\
	$p_1(\alpha_1)\overline{p_1(\alpha_2)}+p_2(\alpha_1)\overline{p_2(\alpha_2)}=p_1(\xi\overline{\alpha_2})\overline{p_1(\alpha_2)}+p_2(\xi \overline{\alpha_2})\overline{p_2(\alpha_2)}.$ \\
	Now using the expressions for $p_1$ and $p_2$ 
		%from Equations (\ref{eq22}) and (\ref{eq23}) 
		we get, 
	\[p_1(\xi\overline{\alpha_2})\overline{p_1(\alpha_2)}+p_2(\xi \overline{\alpha_2})\overline{p_2(\alpha_2)}=\displaystyle \xi a_{11}\overline{\alpha_2}^2+ \xi \overline{a_{12}}\overline{\alpha_2}^3+ \xi^2 a_{12} \overline{\alpha_2}^3+\xi^2  \frac{|a_{12}|^2}{a_{11}} \overline{\alpha_2}^4 + \xi^2 p_{22}^2  \overline{\alpha_2}^4.\] 
	Simplifying, we get
	\[p_1(\xi\overline{\alpha_2})\overline{p_1(\alpha_2)}+p_2(\xi \overline{\alpha_2})\overline{p_2(\alpha_2)}=\xi \overline{\alpha_2}^2(a_{11}+(\overline{a_{12}}+a_{12}\xi )\overline{\alpha_2}+a_{22}\xi \overline{\alpha_2}^2 ).\]

    Let $z_1=a_{11}$, $z_2=(\overline{a_{12}}+a_{12}\xi )\overline{\alpha_2}$ and $z_3=a_{22}\xi \overline{\alpha_2}^2.$ \\
    As $z_1\in \mathbb{R}$, it is enough to prove that $\Im(z_2+z_3)\neq 0.$\\
    Without loss of generality, we assume $\alpha_2$ as that root of $f(z)$, lying outside the closed unit disc, whose argument is smaller. We denote $\arg(\alpha_2)$ by $\theta_2$.
    \begin{enumerate}
        \item [Case (i)]
        Let $\pi/2<\theta<\pi$. We have $z_2=-2b\cos(\theta)(1+\xi)\overline{\alpha_2}$ and hence, $\arg(z_2)=\frac{\theta}{2}-\theta_2$ and $\arg(z_3)=\theta-2\theta_2.$\\
        Recall that $\theta_1=\arg(\alpha_1)$ and we have $\theta_1+\theta_2=\theta.$\\ 
        Further, we have $\arg(1+\xi)=\frac{\theta}{2}$ and $\alpha_1+\alpha_2=a(1+\xi).$ This implies that $\arg(\alpha_1+\alpha_2)=\arg(1+\xi)=\frac{\theta}{2}.$\\
        Since $|\alpha_1|=|\alpha_2|$, $\alpha_1+\alpha_2$ is the angle bisector of the angle between $\alpha_1,\alpha_2$ and thus we must have $0<\arg(\alpha_1)-\arg(\alpha_2)<\pi.$ That is, $0<\theta-2\theta_2<\pi.$  Thus, $0<\frac{\theta}{2}-\theta_2<\frac{\pi}{2}.$ Hence, both $z_2$ and $z_3$ lie in the upper half plane which imply that $\Im(z_2+z_3)\neq 0.$
        
     % Recall that $\theta_1=\arg(\alpha_1)$ and we have $\theta_1+\theta_2=\theta$ and since $\arg(\alpha_2)$ is smaller, we must have $\theta_2<\frac{\theta}{2}.$ Further, since $\alpha_2$ is a continuous function of $\theta$, we must have $\theta_2>0.$ Otherwise, we must have $\theta_2=0$ for some choice of $\theta$ which implies that $\alpha_2\in \mathbb{R}.$ But, this is not true. Hence, we get, $0<\theta_2<\frac{\theta}{2}$ and thus, $\theta>\theta-2\theta_2>0.$ Hence, $z_2,z_3$ lie in upper half plane. Thus, $z_2+z_3\notin \mathbb{R}.$
    \item [Case (ii)] If $t_0<\theta<\pi/2$ then we prove that $\Im(z_2+z_3)\neq 0.$\\
    We have $z_2=-2b\cos(\theta)(1+\xi)\overline{\alpha_2}.$ Note that $\cos\theta>0$ and $\theta_2<\frac{\theta}{2}$ implies that $z_2$ lies in third quadrant. Hence, $\arg(z_2)=\frac{\theta}{2}-\theta_2-\pi.$ \\
    Now, $z_3=a_{22}\xi \overline{\alpha_2}^2$ and thus, $\arg(z_3)=\theta-2\theta_2.$ Hence, $z_3$ lies in first quadrant.\\
    $\Im(z_2)=|z_2|\sin(\frac{\theta}{2}-\theta_2-\pi)$ and $\Im(z_3)=|z_3|\sin(\theta-2\theta_2).$\\
    Suppose $\Im(z_2+z_3)=0$, then we get,\\
    $|z_3|\sin(\theta-2\theta_2)=|z_2|\sin(\frac{\theta}{2}-\theta_2).$
    \\Thus, \begin{equation}\label{eq34}
        \left|\frac{z_3}{z_2}\right|=\frac{1}{2\cos(\frac{\theta}{2}-\theta_2)}<\frac{1}{2\cos(\frac{\theta}{2})}<\frac{1}{\sqrt{2}}.
    \end{equation}
    We have, $|z_2|=2b|1+\xi||\cos(\theta)||\alpha_2|$ and $|z_3|=a_{22}|\alpha_2|^2=a_{22}b.$\\
    Thus, $\left|\frac{z_2}{z_3}\right|=\frac{4\cos(\theta)\cos(\theta/2)|\alpha_2|}{a_{22}}.$\\
    We have, $\cos(\theta)<0.6, \cos(\theta/2)<\sqrt{\frac{1+\cos\theta}{2}}<0.8944$  and \\$ |\alpha_2|=\sqrt{b}<\sqrt{5.46269}=2.3372.$
    \\This gives, $4\cos(\theta)\cos(\theta/2)|\alpha_2|<5.0169.$\\
    We have, $a_{22}=1-3b+2ab-\frac{b}{\Delta}=1+b(2a-3-\frac{b}{\Delta}).$ We know that $5<b<6$. Now $a$ is an increasing function of $\theta$ (see part (2) of Lemma \ref{l6}). This gives $2.5<a<2.5357.$ Further, $0.7<\Delta$ (see part (3) of Lemma \ref{l6}). Thus we get,\\
    $a_{22}=1+b(2a-3-\frac{b}{\Delta})>1+5(2-\frac{1}{0.7})=3.85714.$\\
    Thus, $\left|\frac{z_3}{z_2}\right|>0.7688>\frac{1}{\sqrt{2}}.$ This is a contradiction to (\ref{eq34}).\\
    Hence, we must have $\Im(z_2+z_3)\neq 0.$
    \end{enumerate}
    In the case when $-\pi<\theta<-t_0$, by the symmetry, it can be observed that $\overline{\alpha_2}$ is the root of $f(z)$ corresponding to $-\theta.$ Now, by applying Case (i) and Case(ii) to $-\theta$, we get, $\Im(\overline{z_2+z_3})\neq0$  and thus $\Im(z_2+z_3)\neq0.$
 \\Thus we get, $p_1(\alpha_1)\overline{p_1(\alpha_2)}+p_2(\alpha_1)\overline{p_2(\alpha_2)}\neq 0.$
	Hence, by Corollary \ref{cor1}, the Cauchy dual of the operator $M_z$ is not subnormal.\\
	
{\bf Case 2:} Let $\cos\theta>0.6.$\\
	By Lemma \ref{l3} we know that if $\cos\theta>0.6$ then  $\arg(\alpha_1)=\arg(\alpha_2)$ and $|\alpha_1|\neq |\alpha_2|$. Without loss, suppose $|\alpha_1|<|\alpha_2|$ and $\alpha_2=k\alpha_1$ for some $k>1$. Therefore, note that Corollary \ref{cor1} is not applicable in this case. 
	Here, we apply Theorem \ref{thm3} to reach the desired conclusion. It is enough to prove that for $l=1$, the $2\times2$ minor of the infinite matrix
	\[\displaystyle\sum_{r,t=1}^{2}\left( \frac{1}{a_r\overline{a_t}}\sum_{j=1}^{2}p_j(\alpha_r)\overline{p_j(\alpha_t)}\right)\left(1-\frac{1}{\alpha_r \overline{\alpha_t}}\right)^l\left(\left(\frac{1}{\alpha_r^{m+2}\overline{\alpha_t}^{n+2}}\right)\right)_{m,n\geq 0}\]
	is negative.\\ 
	We have $a_1=\alpha_1-\alpha_2=\alpha_1(1-k)$ and $a_2=\alpha_2-\alpha_1=\alpha_1(k-1)$. For $1\leq r,t\leq 2$, let
	\begin{equation}\label{eq32}
		k_{rt}=\displaystyle\frac{1}{a_r\overline{a_t}}\left(\sum_{j=1}^{2}p_j(\alpha_r)\overline{p_j(\alpha_t)}\right)\left(1-\frac{1}{\alpha_r \overline{\alpha_t}}\right).
	\end{equation}
	For $l=1$, we denote the $2\times2$ principal minor as  
	$M=\begin{vmatrix}
		m_{00}&m_{01}\\
		m_{10}&m_{11}
	\end{vmatrix}$.
	\\We get, 
	\begin{align*}
		m_{00}&=\displaystyle\frac{k_{11}}{|\alpha_1|^4}+\frac{k_{12}}{k^2|\alpha_1|^4}+\frac{k_{21}}{k^2|\alpha_1|^4}+\frac{k_{22}}{k^4|\alpha_1|^4}, \\
		m_{01}&=\displaystyle\frac{k_{11}}{\overline{\alpha_1}|\alpha_1|^4}+\frac{k_{12}}{k^3\overline{\alpha_1}|\alpha_1|^4}+\frac{k_{21}}{k^2\overline{\alpha_1}|\alpha_1|^4}+\frac{k_{22}}{k^5\overline{\alpha_1}|\alpha_1|^4},\\
		m_{10}&=\displaystyle\frac{k_{11}}{\alpha_1|\alpha_1|^4}+\frac{k_{12}}{k^2\alpha_1|\alpha_1|^4}+\frac{k_{21}}{k^3\alpha_1|\alpha_1|^4}+\frac{k_{22}}{k^5\alpha_1|\alpha_1|^4}, \\
		m_{11
		}&=\displaystyle\frac{k_{11}}{|\alpha_1|^6}+\frac{k_{12}}{k^3|\alpha_1|^6}+\frac{k_{21}}{k^3|\alpha_1|^6}+\frac{k_{22}}{k^6|\alpha_1|^6}.
	\end{align*}   
%	Thus, we get, 
%	\begin{align*}
%		M=\displaystyle\frac{1}{|\alpha_1|^{10}}\Bigg[ & k_{11}^2+\displaystyle\frac{k_{11}k_{12}}{k^3}+\frac{k_{11}k_{21}}{k^3}+\frac{k_{11}k_{22}}{k^6}+\frac{k_{12}k_{11}}{k^2}+\frac{k_{12}^2}{k^5}+\frac{k_{12}k_{21}}{k^5}\\
%		&+\frac{k_{12}k_{22}}{k^8}+\frac{k_{21}k_{11}}{k^2}+\frac{k_{21}k_{12}}{k^5}+\frac{k_{21}^2}{k^5}+\frac{k_{21}k_{22}}{k^8}+\frac{k_{22}k_{11}}{k^4}+\frac{k_{22}k_{12}}{k^7}\\
%		&+\frac{k_{22}k_{21}}{k^7}+\frac{k_{22}^2}{k^{10}}-k_{11}^2-\frac{k_{11}k_{12}}{k^2}-\frac{k_{11}k_{21}}{k^3}-\frac{k_{11}k_{22}}{k^5}-\frac{k_{12}k_{11}}{k^3}\\
%		&-\frac{k_{12}^2}{k^5}-\frac{k_{12}k_{21}}{k^6}-\frac{k_{12}k_{22}}{k^8}-\frac{k_{21}k_{11}}{k^2}-\frac{k_{21}k_{12}}{k^4}-\frac{k_{21}^2}{k^5}-\frac{k_{21}k_{22}}{k^7}\\
%		&-\frac{k_{22}k_{11}}{k^5}-\frac{k_{22}k_{12}}{k^7}-\frac{k_{22}k_{21}}{k^8}-\frac{k_{22}^2}{k^{10}}\Bigg]
%	\end{align*}
	Thus, we get, 
	\begin{align*}
		M=\displaystyle\frac{1}{|\alpha_1|^{10}}\Bigg[&\displaystyle\frac{k_{11}k_{22}}{k^6}-\frac{k_{11}k_{22}}{k^5}+\frac{k_{12}k_{21}}{k^5}-\frac{k_{12}k_{21}}{k^6}+\frac{k_{12}k_{21}}{k^5}\\&-\frac{k_{12}k_{21}}{k^4}+\frac{k_{22}k_{11}}{k^4}-\frac{k_{22}k_{11}}{k^5}\Bigg].
	\end{align*}
	Simplifying we get, 
	$M=\displaystyle\frac{\left(1-\frac{1}{k}\right)^2}{k^4|\alpha_1|^{10}}\left(k_{11}k_{22}-k_{12}k_{21}\right).$ \\
	We denote the expression $k_{11}k_{22}-k_{12}k_{21}$ by $G$.\\
	We prove that $G<0$.
	Note that $|a_1|^2=|\alpha_1|^2(k-1)^2$, $|a_2|^2=|\alpha_1|^2(k-1)^2$.
	Also, $a_1\overline{a_2}=-|\alpha_1|^2(k-1)^2=a_2\overline{a_1}$.\\
	We introduce the following notations in order to achieve simplicity in the computations.\\
	Let $x=p_1(\alpha_1)$, $y=p_2(\alpha_1)$ and $w=p_1(\alpha_2)$. Note that $p_2(\alpha_2)=k^2y$.\\
	Substituting the values of $k_{11},k_{12},k_{21},k_{22}$ from equation (\ref{eq32}) we get,
	\begin{align*}
		G=\displaystyle\frac{1}{|\alpha_1|^4(k-1)^4}&\Bigg[\left(|x|^2+|y|^2\right)\left(1-\frac{1}{|\alpha_1|^2}\right)\left(|w|^2+k^4|y|^2\right)\left(1-\frac{1}{k^2|\alpha_1|^2}\right)\\
		&-\left(x\overline{w}+|y|^2k^2\right)\left(1-\frac{1}{k|\alpha_1|^2}\right)\left(w\overline{x}+|y|^2k^2\right)\left(1-\frac{1}{k|\alpha_1|^2}\right)\Bigg].
	\end{align*}
Simplifying we get,
	\begin{align}
		G=\displaystyle\frac{1}{|\alpha_1|^4(k-1)^4}\Bigg[&-\frac{|x|^2|w|^2}{|\alpha_1|^2}\left(1-\frac{1}{k}\right)^2-\frac{k^2|x|^2|y|^2}{|\alpha_1|^2}\left(k^2+1-\frac{1}{|\alpha_1|^2}\right)\nonumber\\
		&-\frac{|y|^2|w|^2}{|\alpha_1|^2}\left(1+\frac{1}{k^2}-\frac{1}{k^2|\alpha_1|^2}\right)-\frac{k^2|y|^4}{|\alpha_1|^2}(k-1)^2\nonumber\\
		&-2|y|^2\Re(x\overline{w})\left(k-\frac{1}{|\alpha_1|^2}\right)^2\Bigg]\label{eq26}\\
		&=\displaystyle\frac{1}{|\alpha_1|^4(k-1)^4}\Bigg[-\frac{|x|^2|w|^2}{|\alpha_1|^2}\left(1-\frac{1}{k}\right)^2-\frac{k^2|y|^4}{|\alpha_1|^2}(k-1)^2-L\Bigg],\nonumber
	\end{align}
	%We now consider the second, third and fifth terms of the above expression with positive sign:
	%We denote the following expression by $L$. 
	where
	$$L=\frac{k^2|x|^2|y|^2}{|\alpha_1|^2}\left(k^2+1-\frac{1}{|\alpha_1|^2}\right)+\frac{|y|^2|w|^2}{|\alpha_1|^2}\left(1+\frac{1}{k^2}-\frac{1}{k^2|\alpha_1|^2}\right)
	+2|y|^2\Re(x\overline{w})\left(k-\frac{1}{|\alpha_1|^2}\right)^2.$$
	
	We now prove that $L>0$. This will in turn prove that $G<0$.\\
	Substituting $y$ in $L$ we have
	\begin{align}
		L=&\frac{k^2|x|^2|p_{22}|^2|\alpha_1|^4}{|\alpha_1|^2}\left(\frac{|\alpha_2|^2+|\alpha_1|^2-1}{|\alpha_1|^2}\right)+\frac{|p_{22}|^2|\alpha_1|^4|w|^2}{|\alpha_1|^2}\left(\frac{|\alpha_2|^2+|\alpha_1|^2-1}{k^2|\alpha_1|^2}\right)\nonumber \\
		&+2|p_{22}|^2|\alpha_1|^4\Re(x\overline{w})\frac{(|\alpha_1||\alpha_2|-1)^2}{|\alpha_1|^4}.\label{eq27}
	\end{align}
Let $\alpha_3=\frac{1}{\overline{\alpha_1}},\alpha_4=\frac{1}{\overline{\alpha_2}}$.  
Since $\alpha_1, \alpha_2, \alpha_3, \alpha_4$ are the roots of $$f(z)=z^4-3(1+\xi)z^3+(8+2\cos\theta)\xi z^2-3(\xi+\xi^2)z+\xi^2,$$
we get,
	$\sum_{i=1}^4\alpha_i=3(1+\xi).$ \\
	Thus, $\sum_{i=1}^4|\alpha_i|\geq \left|\sum_{i=1}^4\alpha_i\right|=|3(1+\xi)|=3\sqrt{2(1+\cos\theta)}$.\\
	But, $\cos\theta>0.6$ gives, $\sum_{i=1}^4|\alpha_i|\geq \displaystyle\frac{12}{\sqrt{5}}$.
	Since, $|\alpha_3|<1,|\alpha_4|<1$, we get, $|\alpha_1|+|\alpha_2|>\displaystyle\frac{12}{\sqrt{5}}-2$.
	Now, $\sqrt{5}<\frac{7}{3}$ and hence, $\frac{12}{\sqrt{5}}-2>3$.
	Thus we get, $|\alpha_1|+|\alpha_2|>3$.
	Also, we have, $(|\alpha_1|-1)(|\alpha_2|-1)>0$ which gives us \[|\alpha_1||\alpha_2|>|\alpha_1|+|\alpha_2|-1>3-1=2.\]
	Thus, we get, $b=|\alpha_1||\alpha_2|>2$ whenever $\cos\theta>0.6$.
	Hence, $|\alpha_1||\alpha_2|-1>1$ which gives us
	\begin{equation}\label{eq28}
		\left(|\alpha_1||\alpha_2|-1\right)^2>|\alpha_1||\alpha_2|-1.
	\end{equation}
	Further, we have, 
	$|\alpha_1|^2+|\alpha_2|^2\geq 2|\alpha_1||\alpha_2|$ and $2|\alpha_1||\alpha_2|\geq |\alpha_1||\alpha_2|$.\\
	Hence, 
	\begin{equation}\label{eq29}
		|\alpha_1|^2+|\alpha_2|^2-1\geq |\alpha_1||\alpha_2|-1.
	\end{equation}
	
	We use inequalities (\ref{eq28}) and (\ref{eq29}) in equation (\ref{eq27}) to get,
	\begin{align*}
		L&\geq\left(|\alpha_1||\alpha_2|-1\right)|p_{22}|^2\left(k^2|x|^2+\frac{|w|^2}{k^2}+2\Re(x\overline{w})\right)\\
		&= \left(|\alpha_1||\alpha_2|-1\right)|p_{22}|^2\left\|kx+\frac{w}{k}\right\|^2>0.
	\end{align*}
	Hence, we get $G<0$, which implies $M<0$. Thus, the Cauchy dual of $M_z$ on $D(\mu)$ is not subnormal.\\
	
	{\bf Case 3:} Suppose $\cos\theta=0.6$.
%	\\We now deal with the case when $\cos\theta=0.6$.\\
	As noted in Part (3) of Lemma \ref{l3}, if  $\cos\theta=0.6$, we have $\alpha_1=\alpha_2=2+i=\alpha(say).$ \\
	We apply Theorem \ref{thm5} to obtain $B=(b_1,b_2)$ such that $D(\mu)=H(B)$ with the equality of norms where $b_j=\frac{p_j}{q}$, for $j=1,2$. By equations (\ref{eq22}) and (\ref{eq23}), we have $p_1(z)=p_{11}z+p_{12}z^2$ and $p_2(z)=p_{22}z^2$, where,\\
	$p_{11}=\sqrt{a_{11}}, p_{12}=\displaystyle\frac{a_{12}}{\sqrt{a_{11}}}$ and $p_{22}=\displaystyle\sqrt{\frac{a_{11}a_{22}-|a_{12}|^2}{a_{11}}}$ with \\
%	\begin{align*}
	$	a_{11}=(a^2-b)(|1+\xi|^2)+\displaystyle\frac{2b}{\Delta},
		a_{12}=-b(1+\xi),
		a_{22}=\displaystyle 1-3b+2ab-\frac{b}{\Delta}. $ \\
%	\end{align*} 
	It can be seen that $b=|\alpha_1\alpha_2|=5,~~~ ~ a=\displaystyle\frac{3b}{b+1}=\frac{5}{2},~~~~|1+\xi|^2=\displaystyle\frac{16}{5}$ and $\Delta=\left[(a-3)\displaystyle\frac{1+\cos\theta}{1-\cos\theta}+\frac{3}{a(1-\cos\theta)}\right]=1.$ \\
	This gives, 
	$a_{11}=14,~~a_{12}=-8-4i,~~~a_{22}=6.$ \\
	At this stage we applied Theorem \ref{thm3} in the Case (2) above. However, here, though we have $D(\mu)=H(B)$, Theorem \ref{thm3} is not applicable as it requires $\alpha_1$ and $\alpha_2$ to be distinct. As stated in \cite[Remark 2.2]{cgr2022}, we need to suitably modify Theorem \ref{thm3} in this case where $B$ has a pole of multiplicity $2$.\\
	 Here, we need to prove that Cauchy dual of $M_z$ on $D(\mu)$ is not subnormal. Therefore, we only need to obtain a condition that will allow us to conclude the same.
	
	In view of \cite[Theorem 3.5]{cgr2022} it is enough to prove that for $l=1$, the $2\times 2$ principal minor of the matrix  
	\begin{equation}\label{eq31}
		\displaystyle\sum_{s=0}^l (-1)^s \binom{l}{s}((B_{m+1+s}B^*_{n+1+s}))_{m,n\geq0}
	\end{equation}
	is negative.\\
	A description of the above matrix and the computations leading to the negativity of the said minor have been carried out below.
	\\We begin by stating a partial fraction formula which can be easily verified:
	\[\displaystyle\frac{a_1+a_2z}{(z-\alpha)^2}=\frac{-a_1/\alpha}{z-\alpha}+\frac{[(a_1+a_2z)/\alpha]z}{(z-\alpha)^2}.\]
	Applying this to polynomials $p_1(z)/z$ and $p_2(z)/z$ we get,
	\[\displaystyle\frac{p_1(z)/z}{(z-\alpha)^2}=\displaystyle\frac{-p_{11}/\alpha}{z-\alpha}+\frac{(p_1(\alpha)/\alpha^2)z}{(z-\alpha)^2} \text{  and  } \displaystyle\frac{p_2(z)/z}{(z-\alpha)^2}=\frac{(p_2(\alpha)/\alpha^2)z}{(z-\alpha)^2}.\]
	Note that 
	$\displaystyle\frac{z}{z-\alpha}=-\displaystyle\sum_{l=0}^\infty \frac{z^{l+1}}{\alpha^{l+1}}\text{  and } \frac{z}{(z-\alpha)^2}=\sum_{l=0}^\infty \frac
	{(l+1)z^{l+1}}{\alpha^{l+2}}.$ \\
	Thus, 
	\begin{align}
		b_1(z)&=\displaystyle\frac{p_1(z)}{q(z)}=\frac{p_1(z)}{(z-\alpha)^2}=\frac{p_{11}}{\alpha}\sum_{l=0}^\infty \frac{z^{l+1}}{\alpha^{l+1}}+\frac{p_1(\alpha)}{\alpha^2}\sum_{l=0}^\infty \frac
		{(l+1)z^{l+1}}{\alpha^{l+2}},\nonumber\\
		b_2(z)&=\displaystyle\frac{p_2(z)}{q(z)}=\frac{p_2(\alpha)}{\alpha^2}\sum_{l=0}^\infty \frac
		{(l+1)z^{l+1}}{\alpha^{l+2}}.\label{eq30}
	\end{align}
	Let $b_j(z)=\displaystyle\sum_{m=1}^\infty b_{j,m}z^m$, $z\in \mathbb{D}$ and $j=1,2$.
	We note here that \\ $B_m=(b_{1,m},b_{2,m})$, $m\geq 1$.\\
	Using equations (\ref{eq30}) we get, 
	$b_{1,m}=\displaystyle\frac{1}{\alpha^{m+1}}\left(p_{11}+\frac{p_1(\alpha)m}{\alpha^2}\right),~~m\geq 1$
	and
	$b_{2,m}=\displaystyle\frac{p_2(\alpha)~m}{\alpha^{m+3}},~~m\geq 1.$ 
	Thus, $\forall m\geq 0, s\geq 0$, we get, \\
	$b_{1,m+1+s}=\displaystyle\frac{1}{\alpha^{m+2+s}}\left(p_{11}+\frac{p_1(\alpha)(m+1+s)}{\alpha^2}\right)$
	and \\
	$b_{2,m+1+s}=\displaystyle\frac{p_2(\alpha)~(m+1+s)}{\alpha^{m+4+s}}.$ 
	After simplification we get, 
	\begin{align*}
		&B_{m+1+s}B^*_{n+1+s}\\
		&=b_{1,m+1+s}\overline{b_{1,n+1+s}}+b_{2,m+1+s}\overline{b_{2,n+1+s}}\\
	%	&=\displaystyle\frac{1}{\alpha^{m+2+s}}\left(p_{11}+\frac{p_1(\alpha)}{\alpha^2}(m+1+s) \right)\frac{1}{\overline{\alpha}^{n+2+s}}\overline{\left(p_{11}+\frac{p_1(\alpha)}{\alpha^2}(n+1+s) \right)}\\
%		&\hspace*{0.5cm}+\displaystyle\frac{p_2(\alpha)~(m+1+s)}{\alpha^{m+4+s}}\times \overline{\frac{p_2(\alpha)~(n+1+s)}{\alpha^{n+4+s}}}\\
%		&=\displaystyle\frac{1}{|\alpha|^{2s}}\frac{1}{\alpha^{m+2}\overline{\alpha}^{n+2}}\Big(|p_{11}|^2+\frac{|p_1(\alpha)|^2}{|\alpha|^4}(m+1+s)(n+1+s)\\
%		&\hspace*{0.5cm}+\displaystyle\frac{p_{11}\overline{p_1(\alpha)}}{\overline{\alpha}^2}(n+1+s)+\frac{\overline{p_{11}}p_1(\alpha)}{\alpha^2}(m+1+s)\Big)\\
%		&\hspace*{0.5cm}+\frac{|p_2(\alpha)|^2(m+1+s)(n+1+s)}{|\alpha|^{2s}\alpha^{m+4}\overline{\alpha}^{n+4}}\\
		&=\displaystyle\frac{|p_{11}|^2}{\alpha^{m+2}\overline{\alpha}^{n+2}}\frac{1}{|\alpha|^{2s}}+\frac{|p_1(\alpha)|^2+|p_2(\alpha)|^2}{\alpha^{m+4}\overline{\alpha}^{n+4}}\Bigg[\frac{(m+1+s)(n+1+s)}{|\alpha|^{2s}} \Bigg]\\
		&\hspace*{0.5cm}+\frac{p_{11}\overline{p_1(\alpha)}}{\alpha^{m+2}\overline{\alpha}^{n+4}}\left(\frac{n+1+s}{|\alpha|^{2s}}\right)+\frac{\overline{p_{11}}p_1(\alpha)}{\alpha^{m+4}\overline{\alpha}^{n+2}}\left(\frac{m+1+s}{|\alpha|^{2s}}\right).
	\end{align*}
%	We will show that, for $l=1$, the $2\times 2$ principal minor of the matrix  
%	\begin{equation}\label{eq31}
%		\displaystyle\sum_{s=0}^l (-1)^s \binom{l}{s}((B_{m+1+s}B^*_{n+1+s}))_{m,n\geq0}
%	\end{equation}
%	is negative and thus by \cite[Theorem 3.5]{cgr2022} we  conclude that the Cauchy dual of $M_z$ on $H(B)$ is not subnormal.\\
	Let the $2\times 2$ principal minor of the matrix in (\ref{eq31}) be denoted by \\$N=\left|(n_{ij})_{0\leq i,j\leq 1}\right|$.
	Hence, 
	\begin{align*}
		n_{00}&= \displaystyle\frac{|p_{11}|^2}{|\alpha|^4}\left(1-\frac{1}{|\alpha|^2}\right)+\frac{|p_1(\alpha)|^2+|p_2(\alpha)|^2}{|\alpha|^8}\left(1-\frac{4}{|\alpha|^2}\right)\\
		&\hspace*{0.5cm}+\frac{p_{11}\overline{p_1(\alpha)}}{|\alpha|^4\overline{\alpha}^2}\left(1-\frac{2}{|\alpha|^2}\right)+\frac{\overline{p_{11}}p_1(\alpha)}{|\alpha|^4\alpha^2}\left(1-\frac{2}{|\alpha|^2}\right),\\
		n_{01}&=\displaystyle\frac{|p_{11}|^2}{|\alpha|^4\overline{\alpha}}\left(1-\frac{1}{|\alpha|^2}\right)+\frac{|p_1(\alpha)|^2+|p_2(\alpha)|^2}{|\alpha|^8\overline{\alpha}}\left(2-\frac{6}{|\alpha|^2}\right)\\
		&\hspace*{0.5cm}+\frac{p_{11}\overline{p_1(\alpha)}}{|\alpha|^4\overline{\alpha}^3}\left(2-\frac{3}{|\alpha|^2}\right)+\frac{\overline{p_{11}}p_1(\alpha)}{|\alpha|^6\alpha}\left(1-\frac{2}{|\alpha|^2}\right),
\end{align*}

\begin{align*}
		n_{10}&=\displaystyle\frac{|p_{11}|^2}{|\alpha|^4\alpha}\left(1-\frac{1}{|\alpha|^2}\right)+\frac{|p_1(\alpha)|^2+|p_2(\alpha)|^2}{|\alpha|^8\alpha}\left(2-\frac{6}{|\alpha|^2}\right)\\
		&\hspace*{0.5cm}+\frac{p_{11}\overline{p_1(\alpha)}}{|\alpha|^6\overline{\alpha}}\left(1-\frac{2}{|\alpha|^2}\right)+\frac{\overline{p_{11}}p_1(\alpha)}{|\alpha|^4\alpha^3}\left(2-\frac{3}{|\alpha|^2}\right),\\
		n_{11}&= \displaystyle\frac{|p_{11}|^2}{|\alpha|^6}\left(1-\frac{1}{|\alpha|^2}\right)+\frac{|p_1(\alpha)|^2+|p_2(\alpha)|^2}{|\alpha|^{10}}\left(4-\frac{9}{|\alpha|^2}\right)\\
		&\hspace*{0.5cm}+\frac{p_{11}\overline{p_1(\alpha)}}{|\alpha|^6\overline{\alpha}^2}\left(2-\frac{3}{|\alpha|^2}\right)+\frac{\overline{p_{11}}p_1(\alpha)}{|\alpha|^6\alpha^2}\left(2-\frac{3}{|\alpha|^2}\right).
	\end{align*}
	Note that
	\[p_{11}=\sqrt{14},~p_{22}=\displaystyle\frac{2}{\sqrt{14}},~~p_1(\alpha)=\frac{20-30i}{\sqrt{14}}~~p_2(\alpha)=\frac{6+8i}{\sqrt{14}}.\]
	Hence, we get,
	\[n_{00}=\displaystyle\frac{228}{625},~~~n_{01}=\frac{344+512i}{3125},~~~n_{10}=\frac{344-512i}{3125},~~~n_{11}=\frac{332}{3125}.\]
	Thus, $N=\displaystyle\frac{-2000}{9765625}<0.$
	Hence, it follows that the Cauchy dual of $M_z$ on $D(\mu)$ is not subnormal.

	Therefore, whenever measure $\mu= \delta_{\zeta_1}+\delta_{\zeta_2}$ where $\zeta_1$ and $\zeta_2$ are two distinct non-antipodal points on $\mathbb{T}$, the Cauchy dual of the operator $M_z$ on corresponding $D(\mu)$ is not subnormal.
\end{proof}

\section{Epilogue}
The authors believe that the main theorem of the present article can be proved for the measure $\mu$ in the general form, i.e. $\mu= c_1\delta_{\zeta_1}+c_2\delta_{\zeta_2},$ where $c_1,c_2$ are positive real numbers, but with the techniques adopted here, the computational complexity increases substantially.    In addition, it would be interesting to extend the main theorem when a measure $\mu$ is supported on three or more points on the unit circle. In particular, a suitable generalization of the notion of points on the unit circle being antipodal needs to be explored.

%The work carried out in this paper confirms the statement conjectured in section 3 of \cite{mkvms2025}. It settles the problem of characterizing the set of all the operators which are cyclic, analytic, $2$-isometries with defect operator of rank $2$ whose Cauchy dual is subnormal. In this context, the following points deserve further exploration:
%\begin{enumerate}
%   \item If the measure $\mu$ is supported at three points $\zeta_1,\zeta_2,\zeta_3$ of $\mathbb{T}$ such that they form an equilateral triangle then Cauchy dual of operator $M_z$ on corresponding $D(\mu)$ is expected to be subnormal. 
%    \item In view of the characterization obtained in this paper, the authors have a feeling that the Cauchy dual of the operator $M_z$ on above $D(\mu)$ is subnormal only in case the support of the measure is at three equidistant points of $\mathbb{T}$.
%\end{enumerate}

\section*{Declarations}

The present work is carried out at the research center at the Department of Mathematics, S. P. College, Pune, India(autonomous).
\section*{Acknowledgments}
The authors are thankful to Dr. Vaidehee Thatte and Dr. Mukul Sholapurkar for fruitful conversations during the preparation of this article. 
\section*{Funding Declaration}
This research has not received any funding.

\printbibliography
\end{document}